\documentclass[11pt]{article}
\usepackage{amsmath}
\usepackage{amsfonts}
\usepackage{graphicx}
\usepackage{latexsym}
\usepackage{url}
\usepackage{color} 
\usepackage{dsfont}
\usepackage{amssymb}
\usepackage{amsthm}
\usepackage{comment}
\usepackage{enumerate}

\excludecomment{en-text}
\includecomment{jp-text}
\excludecomment{comment}

\newcommand{\complexi}{{\tt i}}
\newcommand{\iu}{{\tt i}u}
\newcommand{\iv}{{\tt i}v}

\newcommand{\R}{\mathbb{R}}

\newcommand{\E}{\mathbb{E}}
\newcommand{\PP}{\mathbb{P}}
\newcommand{\N}{\mathbb{N}}

\newcommand{\bbD}{{\mathbb D}}
\newcommand{\bbE}{{\mathbb E}}
\newcommand{\bbF}{{\mathbb F}}
\newcommand{\bbG}{{\mathbb G}}
\newcommand{\bbH}{{\mathbb H}}

\newcommand{\bbN}{{\mathbb N}}

\newcommand{\bbP}{{\mathbb P}}

\newcommand{\bbR}{{\mathbb R}}

\newcommand{\bbZ}{{\mathbb Z}}

\newcommand{\si}{\sigma}

\newcommand{\ep}{\varepsilon}

\newcommand{\De}{\Delta}

\def\nn{\nonumber}

\def\intim{\int_{t_{i-1}}}

\def\sskip{\hspace{0.5cm}}

\newcommand{\toop}{\stackrel{\PP}{\longrightarrow}}
\newcommand{\schw}{\stackrel{d}{\longrightarrow}}
\newcommand{\eqschw}{\stackrel{d}{=}}
\newcommand{\stab}{\stackrel{d_{st}}{\longrightarrow}}

\newcommand{\bee}{\begin{equation}}
\newcommand{\eee}{\end{equation}}
\newcommand{\bea}{\begin{eqnarray}}
\newcommand{\eea}{\end{eqnarray}}
\newcommand{\bean}{\begin{eqnarray*}}
\newcommand{\eean}{\end{eqnarray*}}

\renewcommand{\theequation}{\arabic{section}.\arabic{equation}}

\newtheorem{prop}{Proposition}[section]

\newtheorem{ex}[prop]{Example}
\newtheorem{theo}[prop]{Theorem}
\newtheorem{rem}[prop]{Remark}

\begin{document}

\title{Edgeworth expansion for functionals of continuous diffusion processes}
\author{Mark Podolskij  \thanks{Department of Mathematics, Heidelberg University,  INF 294, 69120 Heidelberg,
Germany, Email: m.podolskij@uni-heidelberg.de.} 
\and Nakahiro Yoshida \thanks{Graduate School of Mathematical
Science, 3-8-1 Komaba, Meguro-ku, Tokyo 153, Japan, Email: nakahiro@ms.u-tokyo.ac.jp}}

\maketitle

\begin{abstract}
This paper presents new results on the Edgeworth expansion for high frequency functionals of continuous
diffusion processes. We derive asymptotic expansions for weighted functionals of the Brownian motion 
and apply them to provide the second order Edgeworth expansion  for power
variation of diffusion processes. Our methodology relies on martingale embedding, Malliavin calculus and stable central limit theorems
for semimartingales.   
Finally, we demonstrate the density expansion for studentized statistics of
power variations.

\ \

\bigskip


\end{abstract}




\section{Introduction} \label{Real-Intro}
\setcounter{equation}{0}
\renewcommand{\theequation}{\thesection.\arabic{equation}}

Edgeworth expansions have been widely investigated by probabilists 
and statisticians in various settings. Nowadays, there exists a vast amount of literature on Edgeworth expansions in the case of independent and identically distributed (i.i.d.) random variables (cf. \cite{H}), weakly dependent variables (cf. \cite{GH})  
or in the framework of martingales (cf. \cite{M, Yoshida1997}). 
We refer to classical books \cite{H,lah} for a comprehensive theory of asymptotic expansions.
We remark that the authors mainly deal with Edgeworth expansions associated with a normal limit.

In the framework of high frequency data (or infill asymptotics), which refers to the sampling scheme in which the time step between two consecutive observations converges to zero while the time span remains fixed, a mixed normal limit appears as a typical 
asymptotic distribution. In the last years a lot of research have been devoted to limit theorems for high frequency observations
of diffusion processes or It\^o semimartingales, see e.g. \cite{BGJPS,  J, J2, KP} among many others. Such limit theorems find
manifold applications in parametric and semiparametric inference for diffusion models, estimation of quadratic variation and related
objects (see e.g. \cite{BS, MA}), testing approaches for semimartingales (see e.g. \cite{AE, DPV}) 
or numerical analysis (see e.g. \cite{JP}).
While asymptotic mixed normality of high frequency functionals has been proved in various settings, the Edgeworth expansions
associated with mixed normal limits have not been considered.

In this paper we present the asymptotic expansion for high frequency statistics of continuous diffusion processes. More 
precisely, we study the second order Edgeworth expansion of weighted functionals of Brownian motion, where the weight 
arises from a continuous SDE, and apply the asymptotic results to power variations of continuous SDE's. Finally, we will 
obtain the density expansion for a studentized version of the power variation. 

Our approach is based on the recent work of Yoshida \cite{Yoshida2012}, who uses a martingale embedding method to obtain
the asymptotic expansion of the characteristic function associated with a mixed normal limit. In a second step 
the asymptotic density expansion is achieved via the Fourier inversion. Let us briefly sketch the main concepts of 
\cite{Yoshida2012}. We are given a functional $Z_n$, which admits the decomposition
\[
Z_n=M_n + r_n N_n,
\]
where $M_n$ is a leading term, 
$r_n$ is a deterministic sequence with $r_n\rightarrow 0$ and $N_n$ is some tight sequence of random variables. Here $M_n$
is a terminal value of a {\it continuous} martingale $(M_t^n)_{t\in [0,1]}$, which converges to a mixed normal limit
in the functional sense. Under various technical conditions, including Malliavin differentiability of the involved objects,
joint stable convergence of $(M_n, N_n)$ and estimates of the tail behaviour of the characteristic function, the paper
\cite{Yoshida2012} demonstrates the second order Edgeworth expansion for the density of $Z_n$ (and, more generally, 
for the density of 
the pair $(Z_n,F_n)$, where $F_n$ is another functional usually used for studentization). The asymptotic theory has been
applied to quadratic functionals $M_n$ in \cite{Yoshida2012, Yoshida2012a}. We would also like to refer to a related
work of \cite{Yoshida1997}, where a martingale expansion in the case of normal limits has been presented.    
It was applied to the Edgeworth expansion for an ergodic diffusion process and an estimator of the
volatility parameter (cf. \cite{DalalyanYoshida2011}). 

Although the paper \cite{Yoshida2012} presents the general theory, its particular application to typical functionals
of continuous diffusion processes is by for not straightforward. When dealing with 
commonly used high frequency statistics, such as e.g. power variations, we are confronted with several levels 
of complications, which we list below: \\
(i) The computation of the second order term $N_n$ in the decomposition of $Z_n$ appears to be rather involved 
(cf. Theorem \ref{th3}). This stochastic second order expansion requires a very precise treatment of the functional $Z_n$. \\
(ii) The joint asymptotic mixed normality of the vector $(M^n, N_n, F_n, C^n)$, where $C^n$ is the quadratic variation
process associated with the martingale $M^n$ and $F_n$ is an external functional mentioned above, is required for
the Edgeworth expansion (cf. Theorem \ref{th6}). The proof of such results relies on stable limit theorems 
for semimartingales (cf. Theorems \ref{generalucp}, \ref{generalclt} and \ref{clt}). \\
(iii) Another ingredients of Edgeworth expansion are the adaptive random symbol $\underline{\sigma}$
and the anticipative random symbol $\overline{\sigma}$ (see \cite{Yoshida2010, Yoshida2012} or Section \ref{sec2}
for the definition of random symbols). While the adaptive random symbol $\underline{\sigma}$ is given explicitly 
using the results of (ii), the anticipative random symbol $\overline{\sigma}$ is defined in an implicit way. 
We will show how this symbol can be determined in Sections \ref{sec32} and \ref{sec33}. For this purpose 
we will apply the Wiener chaos expansion and the duality between the $k$th Malliavin derivative $D^k$ and its
adjoint $\delta^k$. \\
(iv) Checking the technical conditions presented in Section \ref{assumptions} is another difficult task. In particular, 
we need to show the existence of densities and to analyze the tail behaviour of the characteristic function. 
This part involves many elements of Malliavin calculus (cf. Section \ref{sec34} and \ref{sec35}). \\
We see that the derivation of the Edgeworth expansion relies on a combination of various fields of stochastic calculus,
such as limit theorems for semimartingales, Malliavin calculus and martingale methods. We may learn a lot
about the treatment of (iv) from the quadratic case presented in \cite{Yoshida2012a}, but the steps (i)-(iii)
require a completely new treatment in the power variation case. 

The paper is organised as follows. In Section \ref{sec2} we recall the theoretical results of \cite{Yoshida2012}
and demonstrate an application to simple quadratic functionals. Section \ref{sec3} is devoted to functionals 
of Brownian motion with random weights. We will deal with the treatment of the steps (i)-(iv),
although the second order term $N_n$ remains absent. In section \ref{sec4} we show the asymptotic theory 
for the class of generalized power variations of continuous SDE's. In particular, we will determine the asymptotic 
behaviour of the second order term $N_n$. Section \ref{241228-10} combines the results of Sections \ref{sec3} and
\ref{sec4}, and we obtain an Edgeworth expansion for the power variation case. In Section \ref{sec6} we deduce 
the formula for the asymptotic density associated with a studentized version of power variation, which is probably
most useful for applications. Section \ref{proof} is devoted to the derivation of the second order term $N_n$.
Finally, 
Appendix collects the proofs of limit theorems
for semimartingales, which are suitable for functionals considered in this paper.

\section{Asymptotic expansion associated with mixed normal limit} \label{sec2}
\setcounter{equation}{0}
As we are applying various techniques from Malliavin calculus and stable central limit theorems for semimartingales,
we start by introducing some notation. \\ \\
(a) Throughout the paper $\Delta_n$ denotes a sequence of positive real numbers with $\Delta_n\rightarrow 0$ and
such that $1/\Delta_n$ is an integer. 
For the observation times $i\Delta_n$, $i\in \N$, we use a shorthand notation $t_i:=i \Delta_n$. For any function
$f: \R\rightarrow \R$ we denote by $f^{(k)}$ its $k$th derivative; for a function $f:\R^2 \rightarrow \R$ and 
$\alpha =(\alpha_1,\alpha _2)\in \mathbb N^2_{0}$ the operator $d^\alpha$ is defined via 
$d^\alpha= d^{\alpha_1}_{x_1} d^{\alpha_2}_{x_2}$,
where $d^k_{x_i} f$, $i=1,2$, denotes the $k$th partial derivative of $f$. The set $C_p^k(\R)$ (resp. $C_b^k(\R)$)
denotes the space of  $k$ times differentiable functions $f: \R\rightarrow \R$ such that all derivatives up to order
$k$ have polynomial growth (resp. are bounded). Finally, $\complexi:=\sqrt{-1}$.\\ 
(b) The set $\mathbb L^q$ denotes the space of random variables with finite $q$th moment; the corresponding 
$\mathbb L^q$-norms are denoted
by $\|\cdot\|_{\mathbb L^q}$. The notation $Y_n \stab Y$ (resp. $Y_n \toop Y$, $Y_n \schw Y$) stands for stable convergence
(resp. convergence in probability, convergence in law). \\
(c) We now introduce some notions of Malliavin calculus (we refer to the books of Ikeda and Watanabe \cite{IkedaWatanabe1989} 
and Nualart \cite{N} for a detailed exposition of Malliavin calculus). 
Set $\mathbb H=\mathbb L^2([0,1], dx)$ and let 
$\langle \cdot, \cdot \rangle_{\mathbb H}$ denote the usual scalar product on $\mathbb H$. We denote by $D^k$ 
the $k$th Malliavin derivative operator and by $\delta^k$ its unbounded adjoint (also called Skrokhod integral of order $k$). 
The space $\mathbb D_{k,q}$ is the completion of the set of smooth random variables with respect to the norm
\[
\|Y \|_{k,q} := \E[|Y|^q]^{1/q} + \sum_{m=1}^k \E[\|D^m Y\|_{\mathbb H^{\otimes m}}^q]^{1/q}.
\]  
For any $d$-dimensional random variable $Y$ the Malliavin matrix 
is defined via $\si_Y:=(\langle DY_i, DY_j \rangle_{\mathbb H})_{1\leq i,j\leq d}$. We sometimes write $\Delta_Y:=\text{det }\si_Y$
for the determinant of the Malliavin matrix.
Finally, we set $\mathbb D_{k,\infty}= \cap_{q\geq 2} \mathbb D_{k,q}$. 
\\ \\
We start this section by reviewing the theoretical results 
from \cite{Yoshida2010}, which concern the second order Edgeworth expansion 
associated with a mixed normal limit. On a filtered Wiener space $(\Omega, \mathcal F, (\mathcal F_t)_{t\in [0,1]},
\mathbb P)$ we consider a one-dimensional functional $Z_n$, which admits the decomposition
\begin{equation} \label{vndec}
Z_n=M_n + r_n N_n,
\end{equation}
where $r_n$ is a deterministic sequence with $r_n\rightarrow 0$ and $N_n$ is some tight sequence of random variables
(in this paper we will have $r_n=\Delta_n^{1/2}$). 
We assume that the leading term $M_n$ is a terminal 
value of some continuous $(\mathcal F_t)$-martingale $(M_t^n)_{t\in [0,1]}$, i.e. $M_n=M_1^n$. In this paper we are interested in cases,
where $M_n$ (and so $Z_n$) converges stably in law to a mixed normal variable $M$ (stable convergence
has been originally introduced in \cite{REN}). This means:
\begin{equation} \label{mnstab}
M_n \stab M,
\end{equation}
where the random variable $M$ is defined on an extension $(\overline \Omega, \overline{\mathcal F}, 
\overline{\mathbb P})$ of the original probability space $(\Omega, \mathcal F, \mathbb P)$ and, conditionally on $\mathcal F$,
$M$ has a normal law with mean $0$ and conditional variance $C$. In this case we use the notation
\begin{equation*}
M\sim MN(0, C).
\end{equation*}
We recall that a sequence of random variables $(Y_n)_{n\in \N}$ 
defined on $(\Omega, \mathcal F, \mathbb P)$ with values in a metric space $E$ is said 
to converge stably with limit $Y$, written $Y_n \stab Y$, where $Y$ is defined on an extension
$(\overline \Omega, \overline{\mathcal F}, 
\overline{\mathbb P})$ of the original probability space $(\Omega, \mathcal F, \mathbb P)$, 
iff for any bounded, continuous function $g$ and any bounded $\mathcal{F}$-measurable random variable $X$ it holds that
\begin{equation} \label{defstable}
\E[ g(Y_n) X] \rightarrow \overline{\E}[ g(Y) X], \quad n \rightarrow \infty.
\end{equation} 
For statistical applications it is not sufficient to consider the Edgeworth expansion of the law of $Z_n$. It is much more adequate
to study the asymptotic expansion for the pair $(Z_n, F_n)$, where $F_n$ is another functional which converges in probability:
\begin{equation*}
F_n \toop F.
\end{equation*}
When $F_n$ is a consistent estimator of the conditional variance $C$ (i.e. $F=C$), which is the most important application,
we would obtain by the properties of stable convergence:
\begin{equation*}
\frac{Z_n}{\sqrt{F_n}} \schw \mathcal N(0,1).
\end{equation*}
In this case the asymptotic expansion of the law of $(Z_n, F_n)$ would imply the Edgeworth expansion for the studentized 
statistic $Z_n/\sqrt{F_n}$.

We consider the stochastic processes $(M_t)_{t\in [0,1]}$ and $(C^n_t)_{t\in [0,1]}$ with
\begin{equation}
M=M_1, \qquad C_t = \langle M \rangle_t, \qquad C_t^n = \langle M^n \rangle_t, \qquad C_n= \langle M^n \rangle_1.
\end{equation}
Here the process $(M_t)_{t\in [0,1]}$, defined on $(\overline \Omega, \overline{\mathcal F}, 
\overline{\mathbb P})$, represents the stable limit of the continuous $(\mathcal F_t)$-martingale $(M_t^n)_{t\in [0,1]}$,
while $C^n$ denotes the quadratic variation process associated with $M^n$. Now, let us set 
\begin{eqnarray}
\widehat{C}_n= r_n^{-1} (C_n -C), \\[1.5 ex]
\widehat{F}_n= r_n^{-1} (F_n -F).
\end{eqnarray}
Apart from various technical conditions, presented in the Section \ref{assumptions}, our main assumption will be the following:  
\begin{description}
\item[(A1)]
\begin{description}  
\item[(i)] $(M^n_{\cdot}, N_n, \widehat{C}_n, \widehat{F}_n) \stab (M_{\cdot}, N, \widehat{C}, \widehat{F})$.
\item[(ii)] $M_t\sim MN(0, C_t)$.
\end{description} 
\end{description}
\def\bd{\begin{description}}
\def\ed{\end{description}}
\def\dotc{\widehat{C}}
\def\dotf{ \widehat{F}}
\def\dotw{}
\def\calf{\cal{F}}
\newcommand{\cals}{{\cal S}}
\def\bea{\begin{eqnarray}}
\def\eea{\end{eqnarray}}
\def\beas{\begin{eqnarray*}}
\def\eeas{\end{eqnarray*}}
\def\half{\frac{1}{2}}
\def\r{\right}
\def\l{\left}
\def\partialbs{\backslash\!\!\!\partial}
In order to present a second order Edgeworth expansion for the pair $(Z_n, F_n)$ we need to define two {\it random symbols} 
$\underline \sigma$ and $\overline \sigma$, which play a crucial role in what follows. We call $\underline \sigma$ the adaptive
(or classical) random symbol and $\overline \sigma$ the anticipative random symbol.

\subsection{The classical random symbol $\underline \sigma$}
Let $\widetilde{\mathcal F}=\mathcal F \vee \sigma(M)$. We define the random variable
\begin{equation}
\widetilde{C}(M)=\E[\widehat{C}| \widetilde{\mathcal F}].
\end{equation}
In the same way we define the variables $\widetilde{F}(z)$ and $\widetilde N(z)$ such that 
\begin{equation*}
\widetilde{F}(M) = \E[\widehat{F}| \widetilde{\mathcal F}], \qquad \widetilde{N}(M) = \E[N| \widetilde{\mathcal F}].
\end{equation*}
\begin{rem} \label{rem1} \rm
Due to Assumption (A1) (i) we have the pointwise  stable convergence 
$(M_n,  N_n, \widehat{C}_n, \widehat{F}_n) \stab (M, N, \widehat{C}, \widehat{F})$. 
Usually, the limit $(M, N, \widehat{C}, \widehat{F})$ is jointly mixed normal with expectation $\mu \in \R^4$ (and $\mu_1=0$)
and conditional covariance matrix $\Sigma \in \R^{4 \times 4}$. We deduce, for instance, that 
\begin{equation*} 
\widetilde N (M) = \mu_2 + \frac{\Sigma_{12}}{\Sigma_{11}} M.
\end{equation*} 
Consequently, we have $\widetilde{N}(z)=\mu_2 +\frac{\Sigma_{12}}{\Sigma_{11}} z$. 
The quantities $\widetilde{C}(z)$ and $\widetilde F(z)$ are computed similarly. \qed
\end{rem}
Now, the adaptive random symbol $\underline \sigma$ is defined by
\begin{equation} \label{sigmaund} 
\underline \sigma (z, \iu,\iv) = \frac{(\iu)^2}{2} \widetilde{C}(z) +\iu \widetilde N(z) + \iv \widetilde{F}(z).
\end{equation}
Notice that $\underline \sigma$ is a second order polynomial in $(\iu,\iv)$.
The random symbol $\underline \sigma (z,\iu,\iv)$ is called classical, because it appears
already in the martingale expansion in the central limit theorem (\cite{Yoshida1997,Yoshida2001c}), 
i.e. in the case where $C$ is a deterministic constant. In contrast, the anticipative random symbol $\overline \sigma$,
which will be defined in the next subsection, is due to the mixed normality of the limit. In fact, it disappears
if $C$ is non-random.

\subsection{The anticipative random symbol $\overline \sigma$}

The second random symbol $\overline \sigma$ is given in an implicit way. Let 
$\alpha =(\alpha_1,\alpha _2)\in \mathbb N^2_{0}$ with $|\alpha |=\alpha_1+\alpha_2$. Set
\begin{equation*}
\partial^\alpha = \complexi^{-|\alpha|} d^\alpha.
\end{equation*}
We define the quantity $\Phi_n$ by
\begin{equation*}
\Phi_n(u,v) = \E \left[   \exp\Big(- \frac{u^2}{2}C +\iv F \Big) \Big(\mathcal E (\iu M^n)_1 -1 \Big)
\psi_n
\right],
\end{equation*}
where $\mathcal E (H)_t$ denotes the exponential martingale associated with a continuous martingale $H$, i.e.
\begin{equation*}
\mathcal E (H)_t= \exp \Big(H_t -\frac 12 \langle H\rangle_t \Big)= 1+ \int_0^t \mathcal E (H)_s dH_s,
\end{equation*}
and the random variable $\psi_n$ plays a role of a threshold
that ensures the integrability of the above expression, whose precise definition is given in 
Section \ref{assumptions} below.
In particular, $\psi_n$ converges to $1$ in probability.  

\begin{rem} \rm
Recalling the definition of the exponential martingale $\mathcal E (\iu M^n)$, we observe that  $\Phi_n(u,v)$
is closely related to the joint characteristic function of $(M_n, F)$. Condition (A5) of Section \ref{assumptions} 
specifies the tail
behaviour of $\Phi_n(u,v)$. 
When $C=F$ is deterministic, i.e. we are in the framework of a standard central limit theorem, the truncation 
$\psi_n$
can be dropped and we obtain that $\Phi_n(u,v)=0$, since $(\mathcal E (\iu M^n)_t-1)_{t\in [0,1]}$ is a martingale with mean $0$. \qed  
\end{rem} 
Now, we assume that the limit $\Phi^\alpha (u,v):= \lim_{n\rightarrow \infty } r_n^{-1} \partial^\alpha \Phi_n(u,v)$ 
(if it exists)
admits the representation
\begin{equation}
\Phi^\alpha (u,v) = \partial^\alpha \E \left[  \exp\Big(- \frac{u^2}{2}C +\iv F \Big) \overline \sigma(\iu,\iv)  \right], \qquad 
(u,v) \in \R^2,
\end{equation}
where the random symbol $\overline \sigma(\iu,\iv)$ has the form
\begin{equation}
\overline \sigma(\iu,\iv) = \sum_j \overline{c}_j(\iu)^{m_j} (\iv)^{n_j} \qquad \mbox{(finite sum)}
\end{equation}
with $\overline{c}_j\in \cap_{p>1}\mathbb L^p$ (cf. assumption (A4) in Section \ref{assumptions}). We remark that $\overline \sigma(\iu,\iv)$ is a polynomial with random coefficients.

\subsection{ Assumptions and truncation functionals}\label{assumptions}
In this subsection we state the conditions (A2)$_\ell$, (A3),
(A4)$_{\ell,{\mathfrak n}}$, (A5) required in Theorem \ref{th1} below.\footnote{These conditions are the same as or slightly stronger than 
[B2$'$], [B3$'$], [B4] and (12) in \cite{Yoshida2012}.}
Localization techniques will be essential to carry out the computations rigorously. 
We introduce a functional $s_n$ for this purpose. 
\begin{description}
\item[(A2)$_ {\ell}$]
\begin{itemize}
\item[(i)] $F\in\bbD_{\ell+1,\infty}$ and $C\in\bbD_{\ell,\infty}$. 
\item[(ii)] $M_n\in\bbD_{{\ell+1},\infty}$,  $F_n\in\bbD_{\ell+1,\infty}$, $C_n\in\bbD_{\ell,\infty}$, $N_n\in\bbD_{\ell+1,\infty}$  
and $s_n\in\bbD_{\ell,\infty}$.  
Moreover, 
\beas
\sup\Big\{ 
\|M_n\|_{\ell+1,p}+\|\dotc_n\|_{\ell,p}
+\|\dotf_n\|_{\ell+1,p}+\|N_n\|_{\ell+1,p}+\|s_n\|_{\ell,p}
\Big\} <\infty.
\eean
for every $p\geq2$. 
\end{itemize}
\end{description}
\begin{description}
\item[(A3)]
\begin{itemize}
\item[(i)] $ \bbP\l[\Delta_{(M_n,F)}<s_n\r]=O(r_n^{1+\kappa})$ for some positive constant $\kappa$. 
\item[(ii)] For every $p\geq2$, 
\beas 
\limsup_{n\to\infty}\bbE\big[s_n^{-p}\big] <\infty,
\eean
for some $\kappa>0$, 
and moreover $ C^{-1}\in \mathbb{L}^{\infty-}$. Here $\Delta_R$ denotes the determinant of the Malliavin matrix of $R \in \R^d$.
\end{itemize}
\end{description}
\begin{description}
\item[(A4)$_{\ell,{\mathfrak n}}$]
\begin{itemize}
\item[(i)] The random coefficients of the random polynomials 
$\widetilde{C}(z)$, $\widetilde{N}(z)$ and $\widetilde{F}(z)$ are in $\bbD_{4,\infty-}$.
\item[(ii)] The random symbol $\overline{\sigma}$ admits a representation  
\beas
\overline{\sigma}(\iu,\iv)
 =
\sum_j \overline c_j(\iu)^{m_j}(\iv)^{n_j}\hspace{1cm}(\mbox{finite sum}),
\eean
where the numbers $n_j\in \N$ satisfy $n_j\leq{\mathfrak n}$ and $\bar{c}_j\in\bbD_{\ell,\infty-}$.
\end{itemize}
\end{description}
\begin{description}
\item[(A5)] For  some $q\in(1/3,1/2)$, 
\begin{eqnarray*}
\sup_n\sup_{(u,v)\in\Lambda^0_n(2,q)}
|(u,v)|^{3}\Delta_n^{-1/2}|\Phi_n^\alpha(u,v)|&<&\infty
\end{eqnarray*}
for every $\alpha\in\bbZ_+^2$, 
where $\Lambda_n^0(2,q)=\{(u,v)\in\bbR^2;|(u,v)|\leq \Delta_n^{q/2}\}$. 
\end{description}
\vspace*{3mm}

Truncation techniques will play an essential role in derivation of the asymptotic expansion. 
We shall construct a truncation functional $\psi_n$, which has been introduced in the definition
of $\Phi_n(u,v)$, below so that it gives 
uniform convergence of $C^n_t$ and the non-degeneracy of $(Z_n,F_n)$. 
Let $\psi\in C^{\infty}([0,1])$
be a real-valued function with $\psi(x)=1$ for $|x|\leq 1/2$ and $\psi(x)=0$ for $|x|\geq  1$.
Recalling that $C_1=C$, we
define a random variable $\xi_n$ by 
\begin{align} \label{defxi}
\xi_n
&=
10^{-1}\Delta_n^{-c}(C^n_1-C)^2
+2\big[1+4\Delta_{(M^n_1,C)}s_n^{-1}\big]^{-1}+\Delta_n^{c_1}C^2
\\&
+L^*\int_{[0,1]^2}\bigg(\frac{\Delta_n^{-q}|C^n_t-C_t-C^n_s+C_s| }{|t-s|^{3/8}}\bigg)^8dtds, \nonumber
\end{align}
where $L^*$ is a sufficiently large constant, 
$c_1>0$, $c$ satisfies $2q<c<1$ and the constant $q$ is defined in (A5). 
%
%
%
%
Define the $2\times2$ random matrix $R_n'$ by
\beas 
R_n'
&=&
\sigma_{Q_n}^{-1}\bigl(r_n\langle DQ_n,DR_n\rangle_{\mathbb H} 
+ r_n\langle DR_n,DQ_n\rangle_{\mathbb H} +r_n^2\langle DR_n,DR_n\rangle_{\mathbb H}\bigr), 
\eeas
where 
$Q_n=(M_n,F)$ and $R_n=(N_n,\dotf_n)$. 
Obviously 
\bea\label{2109111212-1}
\sigma_{(Z_n,F_n)} &=& \sigma_{Q_n}(I_{2}+R_n'), 
\eea
where $I_2$ is the $2 \times 2$ identity matrix.
Let $\xi_n'= r_n^{-1}|R_n'|^2$. 
We define $\psi_n$ by 
\bea\label{210922-9} 
\psi_n &=& \psi(\xi_n)\psi(\xi_n'). 
\eea

\subsection{The asymptotic expansion of the density of $(Z_n,F_n)$}\label{section.expansion}
We set
\begin{equation} \label{sigmasymb}
\sigma = \underline \sigma + \overline \sigma.
\end{equation}
We remark that due to the definition of $\underline{\sigma}$ and $\overline{\sigma}$ the random symbol $\sigma$
admits the representation 
\begin{equation} \label{sigmarep}
 \sigma(z,\iu,\iv) = \sum_j c_j(z)(\iu)^{m_j} (\iv)^{n_j} \qquad \mbox{(finite sum)}
\end{equation}
for some $c_j(z)\in \cap_{p>1}\mathbb L^p$.  The approximative density of $(Z_n, F_n)$ is defined as 
\begin{eqnarray} \label{pn}
p_n(z,x) &=& \E [ \phi(z;0,C)| F=x] p^F(x) \\[1.5 ex]
&+&r_n \sum_j (-d_z)^{m_j} (-d_x)^{n_j}
\Big( \E \left[  c_j(z)\phi(z;0,C)|F=x\right] p^F(x) \Big) \nonumber, 
\end{eqnarray}
where $p^F$ denotes the density of $F$ and $\phi (\cdot;a,b^2)$ is the density of $\mathcal N(a,b^2)$-distribution. Obviously,  
we will require certain regularity conditions in terms of Malliavin calculus 
in order to validate the existence of the density $p^F$ and the derivatives in (\ref{pn}) as well as to validate the estimate of the approximation error.

For any integrable function $h: \R^2 \rightarrow \R$ we set
\begin{equation} \label{diff}
\Delta_n(h)=\left| \E[h(Z_n,F_n)] - \int h(z,x) p_n(z,x) dzdx \right|.
\end{equation}
The next theorem has been proven in \cite{Yoshida2012}.

\begin{theo} \label{th1} 
Let $\ell=5\vee2[({\mathfrak n}+3)/2]$ with $\mathfrak n = \max_j{n_j}$, where the integers $n_j$ are defined at
\eqref{sigmarep}.  Define the set
$\mathcal{E}(K,\gamma)=\{h:\R^2 \rightarrow \R|~h~\mbox{is measurable and}~ |h(z,x)|\leq K(|z|+|x|)^\gamma\}$ for some 
$K,\gamma >0$. 
Under the assumptions 
(A1), (A2)$_\ell$, (A3), (A4)$_{\ell,{\mathfrak n}}$ and (A5), 
we have that
\begin{equation}
\sup_{h\in \mathcal{E}(K,\gamma)} \Delta_n(h) = o(r_n).
\end{equation}
\end{theo} 
In the following subsection we will explain how this result applies to weighted quadratic  functionals of a Brownian motion.

\subsection{A useful example}
We start by applying the result of Theorem \ref{th1} to a simple example, which however gives a first intuition
how the main quantities are computed. Let $(W_t)_{t\in [0,1]}$ be a 
standard one-dimensional Brownian motion and consider the weighted functional 
\begin{equation} \label{firstmn} 
M_n= \De_n^{-1/2} \sum_{i=1}^{1/\De_n} a(W_{t_{i-1}}) \Big(|\Delta_i^n W|^2 - \Delta_n \Big), \qquad 
\Delta_i^n W=W_{t_{i}} - W_{t_{i-1}},
\end{equation} 
where $a\in C_p^\infty (\R)$ and $\Delta_n\rightarrow 0$ (recall that $1/\De_n$ is an integer). 
In this section we demonstrate how the second order Edgeworth expansion is computed 
for $M_n$. Such quadratic functionals have been already discussed in details 
in \cite{Yoshida2012} and \cite{Yoshida2012a}. For this reason we will dispense with the exact derivation at certain steps of the proof;
in particular, we will not show conditions (A2)$_\ell$-(A5) at this stage. 
We start with the asymptotic properties of the quadratic variation process $C^n$.

\subsubsection{Asymptotic properties of $C^n$}  
The It\^o formula implies the identity  
\begin{equation*}
|\Delta_i^n W|^2 - \Delta_n = 2 \int_{t_{i-1}}^{t_{i}} (W_s-W_{t_{i-1}}) dW_s,
\end{equation*}
and we conclude that $M_n$ is a terminal value of the continuous $(\mathcal F_{t})$-martingale 
\begin{equation} 
M_t^n= 2\De_n^{-1/2} \sum_{i\geq 1} a(W_{t_{i-1}}) 
\int_{t_{i-1}\wedge t}^{t_{i}\wedge t} (W_s-W_{t_{i-1}}) dW_s,
\end{equation}
i.e. $M_n=M_1^n$. We also remark that 
\begin{equation} 
M_t^n= \int_0^t b_s^n dW_s, \qquad b_s^n = 
2\De_n^{-1/2} a(W_{\Delta_n[s/\Delta_n]})  (W_s-W_{\Delta_n[s/\Delta_n]}),
\end{equation}
and thus 
\begin{equation} 
C_t^n = \langle M^n \rangle_t= 4 \De_n^{-1}\int_0^t a^2(W_{\Delta_n[s/\Delta_n]})  (W_s-W_{\Delta_n[s/\Delta_n]})^2 ds.
\end{equation}
Theorem \ref{generalucp} of Appendix implies that 
\begin{equation*} 
C_t^n \toop C_t=2 \int_0^t a^2(W_{s}) ds.
\end{equation*}
\begin{rem} \rm
Notice that Theorems \ref{generalucp} and \ref{generalclt}
of Appendix are formulated for bounded weight functions $a$. 
However, as demonstrated in Section 3 of \cite{BGJPS},
all ingredients involved in our analysis (in particular, the function $a$ in this case) can be assumed to be bounded 
w.l.o.g by a localization technique when proving Theorems \ref{generalucp} and \ref{generalclt}. \qed
\end{rem}
In this example we will consider 
\begin{equation} \label{riemannfn}
F_n= 2\Delta_n \sum_{i=1}^{1/\De_n} a^2(W_{t_{i-1}}),
\end{equation}
which is a Riemann sum approximation of $C=C_1$. We clearly have the convergence in probability  
\begin{equation*} 
F_n \toop C= 2 \int_0^1 a^2(W_{s}) ds.
\end{equation*}
Here and throughout the paper the functional stable convergence  
$(M^n_{\cdot}, N_n, \widehat{C}_n, \widehat{F}_n) \stab (M_{\cdot}, N, \widehat{C}, \widehat{F})$ follows directly
from the general result of Theorem \ref{generalclt} in Appendix, so the assumption (A1) will be always
satisfied. However, for the computation of the random
symbol $\underline{\sigma}$ we only require the pointwise stable convergence 
$(M_n, N_n, \widehat{C}_n, \widehat{F}_n) \stab (M, N, \widehat{C}, \widehat{F})$, which we will present from now 
on. The first convergence of the following proposition is a straightforward consequence of \cite[Section 8]{BGJPS}.

\begin{prop} \label{prop1}
It holds that 
\begin{equation*} 
\Delta_n^{-1/2}(F_n - C) \toop 0.
\end{equation*}
Furthermore, we obtain the stable convergence
\begin{equation*} 
(M_n, \widehat{C}_n) \stab (M,  \widehat{C})\sim MN(0,\Sigma) \qquad \mbox{with} \qquad \Sigma= \int_0^1 \Sigma_s ds, 
\end{equation*}
where the matrix $\Sigma_s$ is defined by 
\begin{equation*} 
\Sigma_s^{11}= 2a^2(W_s), \qquad \Sigma_s^{22}= \frac{16}{3}a^4(W_s), 
\qquad \Sigma_s^{12}=\Sigma_s^{21}=\frac{8}{3}a^3(W_s).
\end{equation*}
\end{prop}

\subsubsection{Computation of $\underline{\sigma}$ and $\overline{\sigma}$}
Now, we start with the computation of the random symbol $\sigma$. We remark that the following
computations are rather typical. The adaptive random symbol $\underline{\sigma}$ is given by
\begin{equation} \label{undersigma}
\underline \sigma (z, \iu,\iv) = \frac{2z(\iu)^2 \int_0^1 a^3(W_s)ds}{3\int_0^1 a^2(W_s)ds} =: z(\iu)^2 \mathcal{C}_1 
\end{equation}
due to Remark \ref{rem1} and $\widetilde N(z)=\widetilde{F}(z)=0$. 

We turn our attention to $\overline{\sigma}$. 
We will not give a rigorous proof, as the details can be found in \cite{Yoshida2012, Yoshida2012a}. Instead we are aiming
to present the most important steps of the derivation (in Section \ref{sec3} we will treat a more general type of functionals
in a detailed manner). Recall that in our case it holds that $F=C$. We set
\begin{equation}
e_t^n(u) = \mathcal E(\iu M^n)_t, \qquad  \Psi(u,v)= \exp\Big((- \frac{u^2}{2} +\iv)C \Big)
\end{equation}
As $e_t^n(u)$ is a continuous exponential martingale we have that
\begin{equation*}
\Phi_n(u,v) = \E \left[   \Psi(u,v) \int_0^1 e_t^n(u) d(\iu M_t^n)  \psi_n\right],
\end{equation*}
for the truncation functional $\xi_n$ defined  
in Section \ref{assumptions}.
The variable $\Delta_n^{-1/2}\Phi_n(u,v)$ has the decomposition: 
\beas 
\Delta_n^{-1/2}\Phi_n(u,v) &=& \check{{\mathfrak A}}_n(u,v)+\hat{{\mathfrak A}}_n(u,v)
\eeas
with 
\begin{align*}
 \check{{\mathfrak A}}_n(u,v)
 &=
\Delta_n^{-1/2} 
\sum_{i=1}^{1/\De_n}
\E \left[   \Psi(u,v) \int_{t_{i-1}}^{t_i} e_{t_{i-1}}^n(u) d(\iu M_t^n)  \psi_n\right], \\[1.5 ex]
\hat{{\mathfrak A}}_n(u,v)
&=
\Delta_n^{-1/2} \sum_{i=1}^{1/\De_n}
 \E \left[   \Psi(u,v)  \int_{t_{i-1}}^{t_i}  (e_t^n(u)-e^n_{t_{i-1}}(u)) d(\iu M_t^n)  \psi_n\right].
\end{align*}
We will see that $\check{{\mathfrak A}}_n(u,v)$ is the dominating term, while $\hat{{\mathfrak A}}_n(u,v)$ turns out to be negligible.  
Setting for simplicity $a_t:=a(W_t)$, we obtain 
\beas 
\check{{\mathfrak A}}_n(u,v)
&=& 
\iu\Delta_n^{-1}\sum_{i=1}^{1/\De_n} \E \left[  \delta^2(1_{I^n_i}^{\otimes2})\times  \Psi(u,v) a_{t_{i-1}}
 \psi_n
e_{t_{i-1}}^n(u) \right],
\eeas
where $I^n_i=1_{(t_{i-1},t_i]}$ and $\delta^2$ denotes the Skorokhod integral of order two. 
Now, we recall the integration by parts (or duality) formula (see e.g. \cite{N}): For any 
$k\in \mathbb N$ and $w\in  \text{Dom}~\delta^k$
and any smooth random variable $Y\in \mathbb D_{k,2}$, it holds that
\begin{align} \label{ibp}
\E[\delta^k(w) Y] = \E[\langle w, D^k Y \rangle_{\mathbb H^{\otimes k}}].
\end{align}  
Applying the duality formula \eqref{ibp} we conclude that  
\beas &&
\check{{\mathfrak A}}_n(u,v)
\\&=& 
\iu\Delta_n^{-1}\sum_{i=1}^{1/\De_n}
\E \left[ 
\int_0^1\int_0^1 1_{I^n_i}^{\otimes2}(s_1,s_2)
D_{s_1,s_2}\big(\Psi(u,v) a_{t_{i-1}}	 \psi_n
e^n_{t_{i-1}}(u)\big)ds_1ds_2 \right]
\\&=& 
\iu\Delta_n^{-1}\sum_{i=1}^{1/\De_n}
\int_{t_{i-1}}^{t_i}\int_{t_{i-1}}^{t_i}ds_1ds_2
\> \E \left[ 
a_{t_{i-1}}e^n_{t_{i-1}}(u)	D_{s_1,s_2}\big(\Psi(u,v)  \psi_n
\big) \right]
\\&=& 
\iu\Delta_n^{-1}\sum_{i=1}^{1/\De_n}
\int_{t_{i-1}}^{t_i}\int_{t_{i-1}}^{t_i}ds_1ds_2
\> \E \left[ 
a_{t_{i-1}}e^n_{t_{i-1}}(u)\psi_n	D_{s_1,s_2}\big(\Psi(u,v) \big) \right]
+o(1),
\eeas
where the term $o(1)$ is explained by the fact that $D\psi_n \rightarrow 0$ since $\psi_n\rightarrow 1$.
Recalling again that $\psi_n \toop 1$, we obtain 
\beas 
\lim_{n\to\infty} \check{{\mathfrak A}}_n(u,v)
&=& 
\int_0^1 \E \left[\iu a_te_t(u)D_{t,t}\Psi(u,v)\right]dt
\eeas
by the functional stable convergence $M^n \stab M$ entailing 
\begin{equation*}
e_{t}^n(u) \stab e_{t}(u) = \exp \Big(\iu M_t +\frac{u^2}{2} C_t \Big). 
\end{equation*}
The process $e_{t}(u)$ is again an exponential martingale and we have 
\begin{equation*}
\overline{\E}[e_{t}(u)| \mathcal F] =1
\end{equation*}
for all $t\geq 0$, $u\in \R$. Hence, 
\beas 
\lim_{n\to\infty} \check{{\mathfrak A}}_n(u,v)
&=& 
\int_0^1 \E \left[\iu a_se_s(u)D_{s,s}\Psi(u,v)\right]ds
\\&=&
\E\big[\Psi(u,v)  \>\iu (4l^2 \mathcal  C_2 + 2l \mathcal C_3)\big]
\eeas
with $l= - \frac{u^2}{2} +\iv$ and
\begin{eqnarray*} 
\mathcal  C_2 &=& \int_0^1 a(W_{s})
  \Big(\int_{s}^1 (a^2)^{\prime } (W_x) dx \Big)^2 ds \\[1.5 ex]
\mathcal  C_3 &=&  \int_0^1 a(W_{s}) \Big(  \int_{s}^1 (a^2)^{\prime \prime} (W_x) dx \Big) ds
\end{eqnarray*}
In a similar way, we obtain the representation 
\beas 
\hat{{\mathfrak A}}_n(u,v)
&=& 
2\iu\Delta_n^{-1}\sum_{i=1}^{1/\De_n}
\int_{t_{i-1}}^{t_i}ds_1 \int_{t_{i-1}}^{s_1} ds_2
\\&&\times
\E \left[ 
D_{s_2}\bigg(D_{s_1}\bigg\{\Psi(u,v) \psi_n a_{t_{i-1}}\times 
(e_{s_1}^n(u)-e^n_{t_{i-1}}(u))\bigg\}\bigg) \right]
\\&\to&
0
\eeas
as $n\to\infty$ for every $(u,v)$ 
by using $\mathbb L^p$-continuity of the densities of the derivatives of $e_{s_1}^n(u)-e^n_{t_{i-1}}(u)$. 
Putting things together we found that 
\beas 
\lim_{n\to\infty} \Delta_n^{-1/2}\Phi_n(u,v)
&=& 
\int_0^1 \E \left[\iu a_se_s(u)D_{s,s}\Psi(u,v)\right]ds
\\&=&
\E\big[\Psi(u,v)  \>\iu (4l^2 \mathcal  C_2 + 2l \mathcal C_3)\big].
\eeas
The corresponding identity for higher order derivatives of $\Phi_n(u,v)$ is proved similarly. Thus,
the anticipative random symbol is given by 
\begin{equation*} 
 \overline \sigma(\iu,\iv)  = \iu (- u^2 +2\iv)^2 \mathcal  C_2 + \iu (- u^2 +2\iv) \mathcal  C_3. 
\end{equation*}
\begin{rem} \label{rem24} \rm
Recall that the statistic $M_n$ depends on $H_2(\Delta_i^n W/ \sqrt{\Delta_n})$, where $H_2$ is the second Hermite
polynomial. In this case the anticipative random symbol $\overline \sigma$ is non-degenerate as we have just proved.  
In Section \ref{sec3} we will show the following fact: Elements of higher order chaos, i.e. $H_m(\Delta_i^n W/ \sqrt{\Delta_n})$
with $m\geq 3$, lead to $\overline \sigma=0$. In other words, for the computation of the anticipative random symbol $\overline \sigma$
of a weighted functional $M_n$ based on $f(\Delta_i^n W/ \sqrt{\Delta_n})$, where $f$ is a measurable function with Hermite
rank at least $2$, only the projection of $f(\Delta_i^n W/ \sqrt{\Delta_n})$ onto the second order Wiener chaos matters. \qed    
\end{rem}
Now, the full random symbol is 
\beas
\sigma (z, \iu,\iv) 
&= &
z(\iu )^2 \mathcal{C}_1 +
  \iu (- u^2 +2\iv)^2 \mathcal  C_2 + \iu (- u^2 +2\iv) \mathcal  C_3. 
\eeas
Therefore the approximative density $p_n(z,x)$ of $(M_n,F_n)$ is given as (recall that $F=C$)
\beas
p_n(z,x) 
&=& \phi(z;0,x) p^C(x) + \Delta_n^{1/2} \Big( 
d_z^2 \{ z\phi(z;0,x) \} p^C(x) \E[\mathcal{C}_1 |C=x]   \\[1.5 ex]
&&-d_z (d_z^2 -2 d_x)^2 \{ \E[\mathcal C_2 \phi(z;0,x) |C=x] p^C(x)\} \\[1.5 ex]
&&-d_z (d_z^2 -2 d_x) \{ \E[\mathcal C_3 \phi(z;0,x) |C=x] p^C(x)\}
\Big).
\eeas

\section{Functionals of Brownian motion with random weights} \label{sec3}
\setcounter{equation}{0}
In this section we go one step further by considering general weighted functionals
of a Brownian motion with weights depending on a given stochastic differential equation, 
and we shall derive an expansion formula. Here the stochastic second order term $N_n$ is still absent. 
In later sections, we will meet an expansion with non-vanishing $N_n$ when considering 
the power variations of diffusion processes. 
However, we will solve two essential problems in this general but concrete situation, 
that is, identification of the anticipative random symbol in this model, and proof of the nondegeneracy of 
the functionals. 
 
On a given Wiener space $(\Omega, \mathcal F, (\mathcal F_t)_{t \in [0,1]}, \mathbb P)$
we consider a $1$-dimensional stochastic differential equation of the form
\begin{equation} \label{sde}
dX_t =  b^{[1]}(X_t) dW_t + b^{[2]}(X_t) dt,  
\end{equation}
where $X_0$ is a bounded random variable, $b^{[1]}, b^{[2]}:  \R \rightarrow \R$ are two deterministic functions
and $W$ is a standard Brownian motion. Sometimes we will use the notation 
\begin{equation*} 
b_t^{[1]}= b^{[1]}(X_t), \qquad b_t^{[2]}= b^{[2]}(X_t).
\end{equation*} 
The somewhat unusual notation $b^{[1]}$, $b^{[2]}$ refers to the fact that the diffusion term $b^{[1]}$ dominates
the drift term $b^{[2]}$ in all asymptotic expansions (so $b^{[1]}$ is the first order term and $b^{[2]}$ is the second order term).
Under standard smoothness conditions the processes $b_t^{[k]}$, $k=1,2$, also satisfy a SDE of the type (\ref{sde}) by It\^o
formula; in this 
case we denote by $b_t^{[k.1]}$ (resp. $b_t^{[k.2]}$) the diffusion term (resp. the drift term) of $b_t^{[k]}$.      
In the same manner we introduce the processes $b_t^{[k_1 \ldots k_d]}$, $k_1,\ldots,k_d=1,2$, recursively.   
We will assume that $b^{[1]}$ and $b^{[2]}$ are in $C^\infty_{b,1}(\bbR)$. 
\footnote{The set of smooth functions such that each derivative of positive order is bounded.}

In this section we consider weighted functionals of the Brownian motion of the type
\begin{equation} \label{weightedw}
M_n= \De_n^{1/2} \sum_{i=1}^{1/\De_n} a(X_{t_{i-1}})f\Big( \frac{\Delta_i^n W}{\sqrt{\Delta_n}}\Big),
\end{equation} 
where $a\in C^\infty_p(\bbR)$ and $f\in C^{11}_p(\bbR)$.
%
Since $f$ has polynomial growth, it holds that $\E[f^2(Z)]<\infty$ with $Z\sim \mathcal N(0,1)$.  
Consequently, the function $f$ exhibits a Hermite expansion. We assume that the function $f$ has the form
\begin{equation} \label{f}
f(x) = \sum_{k=2}^{\infty} \lambda_k H_k(x) \qquad \mbox{in  }\mathbb L^2(\bbR;\phi(x;0,1)dx)  
\end{equation}
with 
\beas
\lambda_k = \frac{\E[f(Z)H_k(Z)]}{k!}, 
\quad Z\sim \mathcal N(0,1),
\eeas
where $H_k$ is the $k$th Hermite polynomial, i.e. $H_0(x)=1$ and 
\begin{equation*}
H_{k}(x)=(-1)^{k}e^{\frac{x^{2}}{2}}{\frac{d^{k}}{dx^{k}}}(e^{-{\frac{%
x^{2}}{2}}}), \qquad k\geq 1.
\end{equation*}
In particular, the Hermite rank of the function $f$ is at least 2 and $\E[f(Z)]=0$ for $Z\sim \mathcal N(0,1)$. 
We will see later that the Hermite rank $1$ would not lead to the asymptotic mixed normal distribution with 
conditional mean $0$.
In this section, we will consider 
\bea\label{241231-10}
F_n= \De_n \mbox{Var}[f(Z)] \sum_{i=1}^{1/\De_n} a^2(X_{t_{i-1}}),
\eea
which is a Riemann sum approximation of 
$C=\langle M\rangle_1$, as the reference variable. It is easy to see the following properties (cf. Proposition \ref{prop1}). 

\begin{prop} \label{prop2}
It holds that 
\begin{equation*} 
F_n \toop C= \mbox{Var}[f(Z)] \int_0^1 a^2(X_s) ds
\end{equation*}
and
\begin{equation*} 
\Delta_n^{-1/2}(F_n - C) \toop 0
\end{equation*}
as $n\to\infty$. 
\end{prop}

\subsection{ A limit theorem for $(M_n,\widehat{C}_n)$ and the adaptive random symbol}
First, we note that for $H= f\Big( \frac{\Delta_i^n W}{\sqrt{\Delta_n}}\Big)$
it holds 
\begin{equation*}
H= \int_0^1 \E[D_s H| \mathcal F_s] dW_s,
\end{equation*}
which is the Clark-Ocone formula.
Consequently, we deduce the identity
\begin{equation*}
f\Big( \frac{\Delta_i^n W}{\sqrt{\Delta_n}}\Big) = \Delta_n^{-1/2}
 \int_{t_{i-1}}^{t_{i}} \E\Big[f'\Big( \frac{\Delta_i^n W}{\sqrt{\Delta_n}}\Big)| \mathcal F_s \Big] dW_s.
\end{equation*}
Thus, we naturally have a continuous square-integrable $({\cal F}_t)$-martingale 
$M^n=(M^n_t)_{t\in[0,1]}$ given by  
\begin{equation} 
M_t^n= \int_0^t b_s^n dW_s, \qquad b_s^n =   a(X_{\Delta_n [s/\Delta_n]}) 
\E\Big[f'\Big( \frac{W_{\Delta_n[s/\Delta_n]+\Delta_n} - W_{\Delta_n[s/\Delta_n]}}{\sqrt{\Delta_n}}\Big)| \mathcal F_s \Big]
\end{equation}
and we deduce that 
\begin{equation} 
C_t^n = \langle M^n \rangle_t=  \int_0^t a^2(X_{\Delta_n[s/\Delta_n]})  
\E^2\Big[f'\Big( \frac{W_{\Delta_n[s/\Delta_n]+\Delta_n} - W_{\Delta_n[ s/\Delta_n]}}{\sqrt{\Delta_n}}\Big)| \mathcal F_s \Big]ds.
\end{equation}
From this identity we obtain the convergence (see Theorem \ref{generalucp} in Appendix) 
\begin{equation*}
C_t^n \toop C_t= \mbox{Var}[f(Z)] \int_0^t a^2(X_s) ds. 
\end{equation*}
By Theorem \ref{generalclt} of Appendix we deduce the following result.

\begin{prop} \label{prop3}
It holds that
\begin{equation*} 
(M_n, \widehat{C}_n) \stab (M,  \widehat{C})\sim MN(0,\Sigma) \qquad \mbox{with} \qquad \Sigma= \int_0^1 \Sigma_s ds, 
\end{equation*}
where the matrix $\Sigma_s$ is defined by
\begin{equation*} 
\Sigma_s^{11}= \mbox{Var}[f(Z)] ~ a^2(X_s), \qquad \Sigma_s^{22}=\Gamma_1 a^4(X_s), \qquad 
\Sigma_s^{12}=\Sigma_s^{21}=\Gamma_2 a^3(X_s),
\end{equation*}
with
\begin{equation*}
\Gamma_1= \mbox{Var} 
\left[ \int_{0}^{1} \E^2[f'(W_1)| \mathcal F_s] ds \right],
\end{equation*} 
\begin{equation*}
\Gamma_2=\text{Cov} 
\left[ f( W_1 ),
\int_{0}^{1} \E^2[f'(W_1)| \mathcal F_s] ds \right].
\end{equation*}
\end{prop}

Notice that the stable convergence in the above proposition does not hold if $f$ has Hermite rank $1$, since in this case
the process $(v_s)_{s\geq 0}$ defined in Theorem \ref{generalclt} is not identically $0$.  
As in the previous section we immediately obtain the adaptive random symbol
\bea\label{241231-2}
\underline \sigma (z, \iu,\iv) 
&=&
 \frac{4z(\iu)^2 \int_0^1 a^3(X_s)ds}{3\mbox{Var}[f(Z)] 
\int_0^1 a^2(X_s)ds} =: z(\iu)^2 \mathcal{C}_1 .
\eea

\subsection{Setting $s_n$}
We  need to define the functionals $s_n$ (and consequently $\xi_n$) to go further.
We set $\beta(x):=Var[(f(Z))] a(x)^2$ with $Z\sim \mathcal N(0,1)$ and $a_t:= a(X_t)$. 
Let 
\beas 
\sigma_{22}(t) 
&=&
\int_0^t \bigg[ \int_r^{1} \beta'_sD_rX_sds\bigg]^2dr. 
\eeas
Define a matrix $\tilde{\sigma}(n,t)$ by 
\beas 
\tilde{\sigma}(n,t)
&=&
\l[\begin{array}{cc}
\tilde{\sigma}_{11}(n,t)&\tilde{\sigma}_{12}(n,t)
\\
\tilde{\sigma}_{12}(n,t)&\sigma_{22}(t)
\end{array}\r]
\eeas
with
\beas 
\tilde{\sigma}_{11}(n,t)&=&
\Delta_n\sum_{i:t_i\leq t}  
\big[
 a_{t_{i-1}} f'(\Delta_n^{-1/2}\Delta^n_iW)\big]^2
\\&&
+\sum_{i:t_i\leq t}  \int_{t_{i-1}}^{t_i}
\bigg[\Delta_n^{1/2}\sum_{k=i+1}^na'_{t_{k-1}}f(\Delta_n^{-1/2}\Delta^n_kW)1_{\{t_k\leq t\}}D_rX_{t_{k-1}}
\bigg]^2dr
\eeas
and 
\beas 
\tilde{\sigma}_{12}(n,t)&=&
\sum_{i:t_i\leq t}  \int_{t_{i-1}}^{t_i} \bigg(
\bigg[
\Delta_n^{1/2}\sum_{k=i+1}^na'_{t_{k-1}}f(\Delta_n^{-1/2}\Delta^n_kW)1_{\{t_k\leq t\}}D_rX_{t_{k-1}}
\bigg] \int_{r}^1 \beta'_sD_rX_sds \bigg)dr
\eeas
for $t\in\Pi^n$. 
Define $s_n$ by 
\beas
s_n &=&
\frac{1}{2}\det \bigg[\tilde{\sigma}\bigg(n,\half\bigg)+\psi\bigg(\frac{m_n}{2{\sf c}_1}\bigg)I_2\bigg],
\eeas
where 
$I_2$ is the $2\times2$ unit matrix, $\psi:\bbR\to[0,1]$ is a smooth function such that 
$\psi(x)=1$ if $|x|\leq1/2$ and $\psi(x)=0$ if $|x|\geq1$, ${\sf c}_1$ is a positive number, and 
\beas 
m_n &=& \Delta_n \sum_{i=1}^{1/\De_n}\big[f'(\Delta_n^{-1/2}\Delta^n_iW)\big]^2. 
\eeas
We will later show that the random variable $s_n$ satisfies assumption (A3).
We define $\xi_n$ using $s_n$ as in Section \ref{assumptions}.

\subsection{Decompositions of the torsion} \label{sec32}

In this subsection we present some preparatory decompositions for the computation of $\overline{\sigma}$. 
Recall that 
$f\in C^{11}_p(\bbR)$ 
and it admits the Hermite expansion $f(x) = \sum_{k=2}^{\infty} \lambda_k H_k(x)$.
Consequently, it holds that $\sum_{k=2}^\infty  k! k^{11} \lambda_k^2 <\infty$.   


The martingale $M^n=(M_t)_{t\in[0,1]}$ admits the local chaos expansion 
\bea\label{241218-1}
M^n_t
&=&  
\Delta_n^{1/2}\sum_{i\geq 1}a_{t_{i-1}}
\sum_{k=2}^\infty k! \lambda_k 
\Delta_n^{-k/2}
\int_{t_{i-1}\wedge t}^{t_i\wedge t}\int_{t_{i-1}}^{s_1}\cdots\int_{t_{i-1}}^{s_{k-1}}
dW_{s_k}\cdots dW_{s_2}dW_{s_1}
\eea
Obviously, each infinite sum in (\ref{241218-1}) is well defined as an $\mathbb L^2$-limit 
when $k\to\infty$.  
Since 
\beas 
H_k(\Delta_n^{-1/2}\Delta^n_iW)
=
k! \Delta_n^{-k/2}
\int_{t_{i-1}}^{t_i}\int_{t_{i-1}}^{s_1}\cdots\int_{t_{i-1}}^{s_{k-1}}dW_{s_k}\cdots dW_{s_2}dW_{s_1}, 
\eeas
we find (\ref{weightedw}) again. 
We recall that for each $(n,u)$, 
the random variable $\sup_{t\in[0,1]}|e^n(u)|$ is bounded uniformly in $\omega$ under 
the truncation by $ \psi_n$. Thus, one can identify $e^n_t(u)$ with 
a stopped $e^n_{\tau_n\wedge t}(u)$ by some stopping time $\tau_n=\tau_n(u)$ that makes 
the stopped process bounded uniformly in $\omega$ for every $n$ (but not uniformly in $n$).  
This remark ensures that variables are in the domain of the Skorokhod integral.
\footnote{A smooth truncation is possible to construct so as to make 
irregularity of the stopping time disappear completely on the remaining event. 
On the other hand, it is also true that one can go without introducing such $\tau_n$ explicitly 
thanks to $\psi_n$  
if the functional $D_{s_2}(e^n_{s_2}D_{s_1}(\Psi(u,v)\psi_n\cdots))$ in the following expressions 
is expanded and interpreted naturally as $e^n_{s_2}\Psi(u,v)\times\cdots$. 
This is always possible because, for every $n$, by some smooth truncation that causes 
$C^n_1\leq A$ locally, the duality operation becomes valid and then the limit $A\to\infty$ 
gives  the formula in expanded form.  }
\begin{en-text}
Define the random variable $\xi_n$ by 
\begin{align} \label{defxi}
\xi_n
&=
10^{-1}\Delta_n^{-c}(C^n_1-C_1)^2
+2\big[1+4\Delta_{(M^n_1,C)}s_n^{-1}\big]^{-1}+\Delta_n^{c_1}C^2
\\&
+L^*\int_{[0,1]^2}\bigg(\frac{\Delta_n^{-q}|C^n_t-C_t-C^n_s+C_s| }{|t-s|^{3/8}}\bigg)^8dtds, \nonumber
\end{align}
where $L^*$ is a sufficiently large constant, 
$c_1>0$ and $c$ satisfies $(2/3<)\>2q<c<1$. Write $\psi_n=\psi(\xi_n)$. 
\end{en-text}

Since the infinite sums in $k$ of (\ref{241218-1}) are also  limits of $\mathbb L^2$-martingales, 
we can validate the exchange of the limit and the sum, and then use the duality between 
the Skorokhod integral and the derivative operator $D$ (cf. \eqref{ibp}) to carry out
\bea\label{241219-1} &&
\sum_{i=1}^{1/\De_n} \E\bigg[\intim^{t_i}e^n_t(u)dM^n_t\>\Psi(u,v)\psi_n\bigg]
\nn\\&=&
\Delta_n^{1/2}\sum_{i=1}^{1/\De_n} \sum_{k=2}^\infty k! \lambda_k\Delta_n^{-k/2}
\E\bigg[\intim^{t_i}e^n_{s_1}(u)\bigg(\intim^{s_1}\cdots\intim^{s_{k-1}}dW_{s_k}\cdots dW_{s_2}\bigg)
dW_{s_1}\>\Psi(u,v)\psi_na_{t_{i-1}}\bigg]
\nn\\&=&
\Delta_n^{1/2}\sum_{i=1}^{1/\De_n} \sum_{k=2}^\infty k! \lambda_k\Delta_n^{-k/2}
\nn\\&&\times
\E\bigg[\int_0^1e^n_{s_1}(u)\bigg(\intim^{s_1}\cdots\intim^{s_{k-1}}dW_{s_k}\cdots dW_{s_2}\bigg)
1_{I^n_i}(s_1)D_{s_1}(\Psi(u,v)\psi_na_{t_{i-1}})ds_1\bigg]
\nn\\&=&
\Delta_n^{1/2}\sum_{i=1}^{1/\De_n} \sum_{k=2}^\infty k! \lambda_k\Delta_n^{-k/2}
\intim^{t_i}ds_1 
\nn\\&&\times
\E\bigg[\intim^{s_1}\bigg(\intim^{s_2}\cdots\intim^{s_{k-1}}dW_{s_k}\cdots dW_{s_3}
 \bigg)dW_{s_2}
e^n_{s_1}(u)D_{s_1}(\Psi(u,v)\psi_na_{t_{i-1}})\bigg]
\nn\\&=&
\Delta_n^{1/2}\sum_{i=1}^{1/\De_n} \sum_{k=2}^\infty k! \lambda_k\Delta_n^{-k/2}
\intim^{t_i}ds_1 \intim^{s_1} ds_2 
\nn\\&&\times 
\E\bigg[\intim^{s_2}\bigg(\intim^{s_3}\cdots\intim^{s_{k-1}}dW_{s_k}\cdots dW_{s_4}\bigg)dW_{s_3} 
D_{s_2}\bigg(e^n_{s_1}(u)D_{s_1}(\Psi(u,v)\psi_na_{t_{i-1}})\bigg)\bigg]. 
\eea
Applying the duality once again, we obtain the decomposition 
\beas 
{\mathfrak A}_n(u,v)&:=&
\Delta_n^{-1/2} \sum_{i=1}^{1/\De_n} \E\bigg[\intim^{t_i}e^n_t(u)dM^n_t\>\Psi(u,v)\psi_n\bigg]
\\&=&
2\sum_{i=1}^{1/\De_n}  \lambda_2\Delta_n^{-1}
\intim^{t_i}ds_1 \intim^{s_1} ds_2 
\E\bigg[
D_{s_2}\bigg(e^n_{s_1}(u)D_{s_1}(\Psi(u,v)\psi_na_{t_{i-1}})\bigg)\bigg]
\\&&
+\sum_{i=1}^{1/\De_n} \sum_{k=3}^\infty k! \lambda_k\Delta_n^{-k/2}
\intim^{t_i}ds_1 \intim^{s_1} ds_2 \intim^{s_2}ds_3
\\&&\times 
\E\bigg[\intim^{s_3}\cdots\intim^{s_{k-1}}dW_{s_k}\cdots dW_{s_4}
\times D_{s_3}\bigg\{
D_{s_2}\bigg(e^n_{s_1}(u)D_{s_1}(\Psi(u,v)\psi_na_{t_{i-1}})\bigg)\bigg\}\bigg]
\\&=&
\ddot{{\mathfrak A}}_n(u,v)+\dddot{{\mathfrak A}}_n(u,v),
\eeas
where 
\begin{align*}
\ddot{{\mathfrak A}}_n(u,v)&=
2\sum_{i=1}^{1/\De_n} \lambda_2\Delta_n^{-1}
\intim^{t_i}ds_1 \intim^{s_1} ds_2 
\E\bigg[
D_{s_2}\bigg(e^n_{s_1}(u)D_{s_1}(\Psi(u,v)\psi_na_{t_{i-1}} )\bigg)\bigg], \\[1.5 ex]
\dddot{{\mathfrak A}}_n(u,v)&=
\sum_{i=1}^{1/\De_n} \Delta_n^{-3/2}
\intim^{t_i}ds_1 \intim^{s_1} ds_2 \intim^{s_2}ds_3 \\[1.5 ex]
&\times 
\E\bigg[
\bigg(\sum_{k=3}^\infty k! \lambda_k\Delta_n^{-(k-3)/2}\intim^{s_3}\cdots\intim^{s_{k-1}}dW_{s_k}\cdots dW_{s_4}\bigg) \\[1.5 ex]
& \times D_{s_3}\bigg\{
D_{s_2}\bigg(e^n_{s_1}(u)D_{s_1}(\Psi(u,v)\psi_na_{t_{i-1}})\bigg)\bigg\}\bigg]. 
\end{align*}
Here we used three times Malliavin differentiability of the objects. We remark that the first term $\ddot{{\mathfrak A}}_n(u,v)$,
which is associated with the second order Wiener chaos, is a dominating quantity, while $\dddot{{\mathfrak A}}_n(u,v)$
will turn out to be negligible (cf. Remark \ref{rem24}).  

\subsection{Identification of the anticipative random symbol} \label{sec33}
We shall specify the limit of ${\mathfrak A}_n(u,v)$. 
First, 
\beas 
|\dddot{{\mathfrak A}}_n(u,v)|
&\leq&
\sum_{i=1}^{1/\De_n} \Delta_n^{-3/2}
\intim^{t_i}ds_1 \intim^{s_1} ds_2 \intim^{s_2}ds_3
\\&&\times 
\bigg\|\sum_{k=3}^\infty k!  \lambda_k\Delta_n^{-(k-3)/2}\intim^{s_3}\cdots\intim^{s_{k-1}}dW_{s_k}\cdots dW_{s_4}
\bigg\|_{\mathbb L^2}
\\&&
\times \bigg\| D_{s_3}\bigg\{
D_{s_2}\bigg(e^n_{s_1}(u)D_{s_1}(\Psi(u,v)\psi_n a_{t_{i-1}})\bigg)\bigg\}\bigg\|_{\mathbb L^2}
\\&\leq&
\Delta_n^{1/2}\sqrt{\sum_{k=3}^\infty k! k^3\lambda_k^2}
\times \sup_{n\in\bbN,\>s_1,s_2,s_3\in[0,1]\atop t_{i-1}<s_3<s_2<s_1\leq t_i}
\bigg\| D_{s_3}\bigg\{
D_{s_2}\bigg(e^n_{s_1}(u)D_{s_1}(\Psi(u,v)\psi_na_{t_{i-1}})\bigg)\bigg\}\bigg\|_{\mathbb L^2}
\\&\to&
0
\eeas
as $n\to\infty$ for every $(u,v)$, since the above supremum is bounded due to assumption (A2$)_8$, and product and
chain rule for the Malliavin derivative. 
Next, we will  treat $\ddot{{\mathfrak A}}_n(u,v)$ essentially  as in Section \ref{sec2} (see also \cite{Yoshida2010}). 
We deform it as 
$\ddot{{\mathfrak A}}_n(u,v)=\check{{\mathfrak A}}_n(u,v)+\hat{{\mathfrak A}}_n(u,v)$ with 
\beas
\check{{\mathfrak A}}_n(u,v)
&=&
2\sum_{i=1}^{1/\De_n} \lambda _2\Delta_n^{-1}\intim^{t_i}ds_1 \intim^{s_1} ds_2 
\E\bigg[a_{t_{i-1}}e^n_{t_{i-1}}(u)D_{s_2}\big(D_{s_1}(\Psi(u,v)\psi_n )\big)\bigg],
\eeas
thanks to $D_sa_{t_{i-1}}=0$ and $D_se^n_{t_{i-1}}(u)=0$ for $s>t_{i-1}$, and 
\beas
\hat{{\mathfrak A}}_n(u,v)
&=&
2\sum_{i=1}^{1/\De_n} \lambda _2\Delta_n^{-1}\intim^{t_i}ds_1 \intim^{s_1} ds_2 
\E\bigg[D_{s_2}\big(e^n_{s_1}(u)-e^n_{t_{i-1}}(u)\big)  \times D_{s_1}(\Psi(u,v)\psi_na_{t_{i-1}} )\bigg]. 
\eeas
Then by continuity of $e^n_\cdot(u)$ in $\bbD_{1,p}$ (see again (A2)), we conclude 
$\hat{{\mathfrak A}}_n(u,v)\to0$ as $n\to\infty$ for every $(u,v)$. 
Since $e^n_{s_1}\Psi(u,v)$ is bounded under truncation by $\psi_n$ or even by its derivative, 
the $\mathbb L^p$-continuity of the objects yields 
\bea\label{241219-5}
\check{{\mathfrak A}}_n(u,v) 
&\to&
\lambda_2 \E\bigg[ \int_0^1   a_t\exp\big(\iu M_t+\half u^2 C_t\big)D_tD_t\Psi(u,v)dt\bigg],
\eea
where
$D_tD_t \Psi(u,v)=\lim_{s\uparrow t}D_sD_t\Psi(u,v)$. 
It should be noted that the integrability and this limiting procedure are valid because 
$D_tD_t \Psi(u,v)=\Psi(u,v) A_t$ with a sum $A_t$ of regular variables, and 
 \beas 
 \mbox{ess sup}_{\omega}\sup_{t\in[0,1]}(C^n_t-C_1)1_{\{|\xi_n|<1\}}\leq \Delta_n^{c/2}\leq1<\infty
\eeas 
for all $n$, due to $C^n_t\leq C^n_1$ and the construction of the quantity $\xi_n$ in Section \ref{assumptions}. 
Furthermore, 
\beas &&
\lambda_2 \E\bigg[\int_0^1   a_t\exp\big(\iu M_t+\half u^2 C_t\big)D_tD_t\Psi(u,v)dt\bigg]
\\&=&
\lambda_2 \E\bigg[ \int_0^1  a_t
\E\big[\exp\big(\iu M_t\big)|{\cal F}\big] \big\{\exp\big(\half u^2 C_t\big)\times \Psi(u,v)\big\}A_t dt \bigg]
\\&=&
\lambda_2 \E\bigg[ \int_0^1  a_tD_tD_t\Psi(u,v)dt\bigg].
\eeas
Consequently, for 
\beas 
\Phi^{\alpha}_n(u,v)
&=&
\complexi^{-|\alpha|}d_{(u,v)}^\alpha
\E\big[L^n_1(u)\Psi(u,v)\psi_n\big]
=
\E\bigg[\iu\int_0^1e^n_t(u)dM^n_t\>\Psi(u,v)\psi_n\bigg],
\eeas
where $L^n_t(u)=e^n_t(u)-1$, we obtain
\beas 
\tilde{\Phi}^{\alpha}(u,v)
&=&
\lim_{n\to\infty}\Delta_n^{-1/2}\Phi^{\alpha}_n(u,v)
\\&=&
\lim_{n\to\infty}\complexi^{-|\alpha|}d_{(u,v)}^\alpha
\big(\iu{\mathfrak A}_n(u,v)\big)
\\&=&
\lambda_2 \complexi^{-|\alpha|}d_{(u,v)}^\alpha
\E\bigg[ \int_0^1 \iu a_tD_tD_t\Psi(u,v)dt\bigg]
\\&=&
\lambda_2 \E\bigg[\Psi(u,v)\cdot  \int_0^1 \iu a_t
\bigg(
(-\frac{u^2}{2}+\iu)^2(D_tC)^2
+(-\frac{u^2}{2}+\iu)D_tD_tC
\bigg)dt\bigg]
\eeas
Therefore, 
\bea\label{241231-1} 
\overline{\sigma}(\iu,\iv)
&=&
\lambda_2 \int_0^1 \iu a_t
\bigg(
(-\frac{u^2}{2}+\iu)^2(D_tC)^2
+(-\frac{u^2}{2}+\iu)D_tD_tC
\bigg)dt.
\eea
We recall that the process $DX_t$ is given as the solution of the SDE
\beas
D_s X_t&=& b^{[1]}(X_s) + \int_s^t (b^{[2]})'(X_{u}) D_s X_{u} du + \int_s^t (b^{[1]})'(X_{u}) D_s X_{u} dW_{u}
\eeas
for $s\leq t$ (and $0$ when $s>t$), and 
\beas 
D_rD_sX_t
&=&
(b^{[1]})'(X_s)D_rX_s+\int_x^t(b^{[2]})''(X_u)D_rX_{u}D_sX_{u}du
+\int_s^t(b^{[2]})'(X_{u})D_rD_sX_{u}du
\\&&
+\int_s^t(b^{[1]})''(X_{u})D_rX_{u}D_sX_{u}dW_{u}+\int_s^t(b^{[1]})'(X_{u})D_rD_sX_{u}dW_{u}
\eeas
for $r<s\leq t$. 
Then (\ref{241231-1}) implies the identity.
\begin{en-text}
As before we obtain that  
\begin{eqnarray*} 
\Phi_n(u,v) = \Delta_n^{1/2} \E[ \iu\lambda_2 \Psi(u,v) (l^2 \mathcal  C_2 + l \mathcal C_3+ l \mathcal C_4)
] + o(\Delta_n^{1/2}),
\end{eqnarray*} 
with 
\end{en-text}
\bea\label{241231-3}
\overline \sigma (\iu,\iv)
\nn&=&
 \iu\lambda_2  \Big((- \frac{u^2}{2} +\iv)^2\mbox{Var}^2[f(Z)] 
\mathcal  C_2 + (- \frac{u^2}{2} +\iv)\mbox{Var}[f(Z)] (\mathcal C_3+  \mathcal C_4) \Big)
\nn\\&&
\eea
with 
\beas
\mathcal  C_2 &=& \int_0^1 a(X_{s})
  \Big(\int_{s}^1 (a^2)^{\prime } (X_u) D_{s} X_u du \Big)^2 ds \\[1.5 ex]
\mathcal  C_3 &=&  \int_0^1 a(X_{s}) \Big( \int_{s}^1 (a^2)^{\prime \prime} (X_u) (D_s X_u)^2 du  \Big) ds \\[1.5 ex]
\mathcal  C_4 &=&  \int_0^1 a(X_{s}) \Big( \int_{s}^1 (a^2)^{\prime } (X_u) D_s D_s X_u du  \Big) ds.
\eeas
Now, having obtained the full random symbol $\sigma=\underline{\sigma}+\overline{\sigma}$ 
and hence the density $p_n(z,x)$ for $\sigma$, we can formulate 
the following statement, which generalizes the results of \cite[Theorem 6]{Yoshida2012} and  \cite[Theorem 1]{Yoshida2012a} 
on the quadratic form to the weighted power variation of Brownian motion.

\begin{theo}\label{241227-8} 
Let $b^{[1]},b^{[2]}\in C^\infty_{b,1}(\bbR)$,  
$a\in C^\infty_p(\bbR)$ and 
$f\in C^{11}_p(\bbR)$.
Let the functional $F_n$ be given by (\ref{241231-10}). 
Define $\beta(x)=\mbox{Var}[f(Z)]\>a(x)^2$ for a standard normal random variable $Z$. 
Assume that the following conditions are satisfied: 
\begin{description} 
\item[(C1)] 
$\displaystyle\inf_{x} |b^{[1]}(x)|>0$ 
 and $\displaystyle\inf_{x} |a(x)|>0$. 
\item[(C2)] 
$\sum_{k=1}^{\infty}|\beta^{(k)}(X_0)|>0$.
\end{description}
Then for any positive numbers $K$ and $\gamma$, 
it holds that 
\beas 
\sup_{h\in\mathcal{E}(K,\gamma)}
\big|\E[f(M_n,F_n)]-\int h(z,x)p_n(z,x)dzdx\big|
&=&
o(\sqrt{\Delta_n}\>)
\eeas
as $n\to\infty$, where the set $\mathcal{E}(K,\gamma)$ was defined in Theorem \ref{th1}.    
\end{theo}

In the rest of this section, we will prove Theorem \ref{241227-8}. 
We will  verify conditions (A1), (A2)$_\ell$, (A3), (A4)$_{\ell,{\mathfrak n}}$ and (A5) of 
Theorem \ref{th1} for $\ell=10$. The conditions of Theorem \ref{241227-8} trivially imply 
(A1) and  (A2)$_\ell$. We already have (A4)$_{\ell,{\mathfrak n}}$. In the following
subsections we concentrate on proving (A3) and (A5).

\subsection{Estimate of the characteristic functions} \label{sec34}

We shall now show condition (A5) of Section \ref{assumptions} 
under the assumptions of Theorem \ref{241227-8}, namely 
%
\bea\label{241216-1}
\sup_n\sup_{(u,v)\in\mathcal{E}^0_n(2,q)}
|(u,v)|^{3}\Delta_n^{-1/2}|\Phi_n^\alpha(u,v)|&<&\infty
\eea
for 
\beas 
\Phi_n^\alpha(u,v) 
&=& 
\complexi^{-|\alpha|}d_{(u,v)}^\alpha \E[L^n_1(u)\Psi(u,v)\psi_n], \qquad L^n_t(u)=e^n_t(u)-1.
\eeas
We shall take a similar way as that was done in the proof of   
\cite[Theorem 4]{Yoshida2010}.  
Though only the quadratic form was treated there, the same method works in our situation. 
The idea is to apply the duality formula twice and use nondegeneracy of 
the Malliavin matrix of $(M^n_t,F)$ together with 
that of $C-C_t$, in the expression 
\beas 
\Phi_n^\alpha(u,v) 
&=&
\complexi^{-|\alpha|}d_{(u,v)}^\alpha \E\bigg[\int_0^1e^n_t(u)d(\iu M^n_t)\>\Psi(u,v)\psi_n\bigg].
\eeas
For this purpose, the representation (\ref{241219-1}) is useful.  
By the $\mathbb L^2$-convergence, we see that
\beas
&&
\sum_{i=1}^{1/\De_n} \E\bigg[\intim^{t_i}e^n_t(u)dM^n_t\>\Psi(u,v)\psi_n\bigg]
\nn\\&=&
\Delta_n^{1/2}\sum_{i=1}^{1/\De_n} \sum_{k=2}^\infty k! \lambda_k\Delta_n^{-k/2}
\intim^{t_i}ds_1 \intim^{s_1} ds_2 
\nn\\&&\times 
\E\bigg[\intim^{s_2}\intim^{s_3}\cdots\intim^{s_{k-1}}dW_{s_k}\cdots dW_{s_4}dW_{s_3} 
\>a_{t_{i-1}}D_{s_2}\bigg(e^n_{s_1}(u)D_{s_1}(\Psi(u,v)\psi_n)\bigg)\bigg]
\nn\\&=&
\Delta_n^{-1/2}\sum_{i=1}^{1/\De_n} 
\intim^{t_i}ds_1 \intim^{s_1} ds_2 
\nn\\&&\times 
\E\bigg[\sum_{k=2}^\infty k! \lambda_k\Delta_n^{-(k-2)/2}
\intim^{s_2}\intim^{s_3}\cdots\intim^{s_{k-1}}dW_{s_k}\cdots dW_{s_4}dW_{s_3} 
\>a_{t_{i-1}}D_{s_2}\bigg(e^n_{s_1}(u)D_{s_1}(\Psi(u,v)\psi_n)\bigg)\bigg]
\nn\\&=&
\Delta_n^{-1/2}\sum_{i=1}^{1/\De_n} 
\intim^{t_i}ds_1 \intim^{s_1} ds_2 
\>\E\bigg[f^\dagger_{n,i,s_2}
\>a_{t_{i-1}}D_{s_2}\bigg(e^n_{s_1}(u)D_{s_1}(\Psi(u,v)\psi_n)\bigg)\bigg],
\eeas
where 
\beas 
f^\dagger_{n,i,s_2}
&=&
\sum_{k=2}^\infty k! \lambda_k\Delta_n^{-(k-2)/2}
\intim^{s_2}\intim^{s_3}\cdots\intim^{s_{k-1}}dW_{s_k}\cdots dW_{s_4}dW_{s_3} , 
\eeas
and consequently reach the representation
\bea\label{241220-2}
\iu \sum_{i=1}^{1/\De_n} \E\bigg[\intim^{t_i}e^n_t(u)dM^n_t\>\Psi(u,v)\psi_n\bigg]
&=&
\Delta_n^{-1/2}\sum_{i=1}^{1/\De_n} 
\intim^{t_i}ds_1 \intim^{s_1} ds_2 
\>E^n_i(u,v)_{s_1,s_2}
\eea
where
\bea\label{241220-5}
E^n_i(u,v)_{s_1,s_2}
&=&
iu \>\E\bigg[f^\dagger_{n,i,s_2}
\>a_{t_{i-1}}D_{s_2}\bigg(e^n_{s_1}(u)D_{s_1}(\Psi(u,v)\psi_n)\bigg)\bigg].
\eea
Let 
\beas 
\bbE^n_s(u,v)
&=&
e^n_s(u)\Psi(u,v).
\eeas
Then $\bbE_s(u,v)$ has the {FGH-decomposition} (cf. \cite[page 22]{Yoshida2012}): 
\beas 
\bbE^n_s(u,v)
&=&
\bbF^n_s(u,v)\bbG_s(u)\bbH^n_s(u)
\eeas
with
\beas 
\bbF^n_s(u,v) &=& \exp\big(\iu M^n_s+\iv C_1\big),\sskip C_1=C,
\\
\bbG_s(u) &=& \exp\bigg(-\half u^2(C_1-C_s)\bigg),
\\
\bbH^n_s(u) &=& \exp\bigg(\half u^2(C^n_s-C_s)\bigg).
\eeas
{From} (\ref{241220-5}) and the FGH-decomposition, 
\bea\label{241220-6}
E^n_i(u,v)_{s_1,s_2}
&=&
\E\bigg[\bbF^n_{s_1}(u,v)\bbG_{s_1}(u)\bbH^n_{s_1}(u)\psi^n_{s_1,s_2}(u,v)
f^\dagger_{n,i}\>a_{t_{i-1}}\bigg],
\eea
where 
\beas 
\psi^n_{s_1,s_2}(u,v)
&=&
\iu \>\big(e^n_{s_1}(u)\Psi(u,v)\big)^{-1}\> D_{s_2}\big(e^n_{s_1}(u)D_{s_1}(\Psi(u,v)\psi_n)\big)
\\&=&
\bigg\{\psi_n\big(-\frac{u^2}{2}+\iv\big)D_{s_1}C_1+D_{s_1}\psi_n\bigg\}
\iu\bigg(\iu D_{s_2}M^n_{s_1}+\frac{u^2}{2}D_{s_2}C^n_{s_1}\bigg)
\\&&
+\psi_n  \iu \big(-\frac{u^2}{2}+\iv\big)^2(D_{s_2}C_1)(D_{s_1}C_1)
\\&&
+2(D_{s_2}\psi_n)\iu\big(-\frac{u^2}{2}+\iv\big)D_{s_1}C_1
+D_{s_2}D_{s_1}\psi_n \>\iu.
\eeas
Suppose that the following condition, which we will prove in the next subsection, is satisfied
for $\ell=10$: 
\begin{description}
\item[(C2$^\flat$)] 
The variables $s_n$ ($n\in\bbN$) satisfy the following conditions.
\begin{description}
\item[(i)] 
$\displaystyle
\sup_{t\geq\half} \mathbb P\big[\det \sigma_{(M^n_t,C_1)}<s_n\big]
=O(\Delta_n^{4/3+\ep})$ as $n\to\infty$ for some $\ep>0$. 
\item[(ii)] 
$\displaystyle 
\limsup_{n\to\infty}\E[s_n^{-p}]<\infty$ for every $p>1$. 
\item[(iii)] 
$\displaystyle \limsup_{n\to\infty}\|s_n\|_{\ell,p}<\infty$ for every $p\geq2$. 
\end{description}
\end{description}
Note that condition (C2$^\flat$) immediately implies (A3).  
Now following the (a)-(h) procedure of \cite[page 22]{Yoshida2012}
and the argument of the proof of Theorem 4 therein, we can obtain 
\bea\label{241220-10}
\sup_n\sup_{i=1,...,n}\sup_{s_1,s_2:t_{i-1}<s_1<s_2\leq t_i} \sup_{(u,v)\in\Lambda^0_n(2,q)}
|(u,v)|^{3}\big| E^n_i(u,v)_{s_1,s_2}\big|
&<& \infty
\eea
by applying the integration-by-parts formula at most $8$ 
times.  
More precisely, we introducing a new truncation 
\beas 
\psi_{n,s_1}
&=& 
\psi\bigg(2\big[1+4\Delta_{(M^n_{s_1},C)}s_n^{-1}\big]^{-1}\bigg),
\eeas
which will be used when the integration-by-parts formula for $(M^n_{s_1},C)$ is applied 
for $s_1\geq 1/2$. 
We have the decomposition of $E^n_i(u,v)_{s_1,s_2}$ expressed by (\ref{241220-6}): 
\beas &&
E^n_i(u,v)_{s_1,s_2}
\\&=&
\E\bigg[\bbF^n_{s_1}(u,v)\bbG_{s_1}(u)\bbH^n_{s_1}(u)\psi^n_{s_1,s_2}(u,v) \psi^n_{s_1}
f^\dagger_{n,i}\>a_{t_{i-1}}\bigg]
+R_{n,s_1,s_2}(u,v)
\eeas
with 
\beas 
|R_{n,s_1,s_2}(u,v)| 
&\leq&
K\Delta_n^{-5q/2}\sup_{s'}\|1-\psi^n_{s'}\|_{\mathbb L^p}
\eeas
for all $n,s_1$ and restricted $(u,v)$, where $K$ a constant. 
The right-hand side can be shown to be  of order $o(\Delta_{n}^{3q/2})$ for sufficiently small numbers $q>1/3$ (cf.
assumption (A5)) and $p>1$. 
Then, as already noticed, we can follow the (a)-(h) procedure of \cite{Yoshida2010}, 
by using the FGH-decomposition, but with $\psi(\xi_n)\psi_{n,s_1}$ for truncation, 
to obtain (\ref{241220-10}).

Finally, we obtain (\ref{241216-1}) for $\alpha=0$ from (\ref{241220-10}). 
When $\alpha\not=0$, the argument of the proof is essentially the same as above. 
As a conclusion, (\ref{241216-1}) (and consequently (A5)) holds for every $\alpha$ under the assumptions $(C1)$ and $(C2^\flat)$. 

Obviously, condition (A3) is valid under $(C1)$ and $(C2^\flat)$. 
In particular, the non-degeneracy of $C$ simply follows from $\inf_x|a(x)|>0$. Thus, we are left to proving condition 
$(C2^\flat)$.

\subsection{Proof of (C2$^\flat$)} \label{sec35}
We shall now prove that condition (C2$^\flat$) holds under the assumptions of Theorem \ref{241227-8}. 
Recall that 
\beas 
M^n_t &=& 
\Delta_n^{1/2}\sum_{i:t_i\leq t} a(X_{t_{i-1}})f(\Delta_n^{-1/2}\Delta^n_iW)
\eeas
for $t\in\Pi^n=\{t_i\}$. 
We deduce that
\beas 
D_rM^n_t 
&=& 
\sum_{i:t_i\leq t}  a_{t_{i-1}} f'(\Delta_n^{-1/2}\Delta^n_iW)1_{(t_{i-1},t_i]}(r)
\\&&
+
\Delta_n^{1/2}\sum_{i:t_i\leq t}  a'_{t_{i-1}}D_rX_{t_{i-1}}f(\Delta_n^{-1/2}\Delta^n_iW)1_{\{r\leq t_{i-1}\}}
\\&=&
\sum_{i:t_i\leq t}  
\bigg[
 a_{t_{i-1}} f'(\Delta_n^{-1/2}\Delta^n_iW)
\\&&
+\Delta_n^{1/2}\sum_{k=i+1}^na'_{t_{k-1}}f(\Delta_n^{-1/2}\Delta^n_kW)1_{\{t_k\leq t\}}D_rX_{t_{k-1}}
\bigg]
1_{(t_{i-1},t_i]}(r)
\eeas
for $t\in\Pi^n$, 
where $\sum_{k=n+1}^n\cdots=0$. 
Hence
\beas 
\sigma_{11}(n,t) &:=&\sigma_{M^n_t}
\\&=&
\sum_{i:t_i\leq t}  \int_{t_{i-1}}^{t_i}
\bigg[
 a_{t_{i-1}} f'(\Delta_n^{-1/2}\Delta^n_iW)
\\&&
+\Delta_n^{1/2}\sum_{k=i+1}^na'_{t_{k-1}}f(\Delta_n^{-1/2}\Delta^n_kW)1_{\{t_k\leq t\}}D_rX_{t_{k-1}}
\bigg]^2dr
\eeas
for $t\in\Pi^n$. 
%
We have 
$C_t =\int_0^t \beta(X_s)ds$. 
Since 
\beas 
D_r C_t &=& \int_r^t \beta'_sD_rX_sds,\sskip t\in[0,1],  
\eeas
we obtain
\beas 
\sigma_{12}(n,t) 
&:=&
\langle DM^n, DC\rangle_{\mathbb H}
\\&=&
\sum_{i:t_i\leq t}  \int_{t_{i-1}}^{t_i} \bigg(
\bigg[
 a_{t_{i-1}} f'(\Delta_n^{-1/2}\Delta^n_iW)
\\&&
+\Delta_n^{1/2}\sum_{k=i+1}^na'_{t_{k-1}}f(\Delta_n^{-1/2}\Delta^n_kW)1_{\{t_k\leq t\}}D_rX_{t_{k-1}}
\bigg] \int_{r}^1 \beta'_sD_rX_sds \bigg) dr
\eeas
for $t\in\Pi^n$.  
The Malliavin  matrix of $(M^n_t,C)$ is 
\beas 
\sigma_{(M^n_t,C)}
&=&
\l[\begin{array}{cc}
\sigma_{11}(n,t)&\sigma_{12}(n,t)
\\
\sigma_{12}(n,t)&\sigma_{22}(1)
\end{array}\r]
\eeas
for $t\in\Pi^n$. 
Let 
\beas 
\sigma(n,t)
&=&
\l[\begin{array}{cc}
\sigma_{11}(n,t)&\sigma_{12}(n,t)
\\
\sigma_{12}(n,t)&\sigma_{22}(t)
\end{array}\r].
\eeas
\begin{en-text}
Define a matrix $\tilde{\sigma}(n,t)$ by 
\beas 
\tilde{\sigma}(n,t)
&=&
\l[\begin{array}{cc}
\tilde{\sigma}_{11}(n,t)&\tilde{\sigma}_{12}(n,t)
\\
\tilde{\sigma}_{12}(n,t)&\sigma_{22}(t)
\end{array}\r]
\eeas
with
\beas 
\tilde{\sigma}_{11}(n,t)&=&
\Delta_n\sum_{i:t_i\leq t}  
\big[
 a_{t_{i-1}} f'(\Delta_n^{-1/2}\Delta^n_iW)\big]^2
\\&&
+\sum_{i:t_i\leq t}  \int_{t_{i-1}}^{t_i}
\bigg[\Delta_n^{1/2}\sum_{k=i+1}^na'_{t_{k-1}}f(\Delta_n^{-1/2}\Delta^n_kW)1_{\{t_k\leq t\}}D_rX_{t_{k-1}}
\bigg]^2dr
\eeas
and 
\beas 
\tilde{\sigma}_{12}(n,t)&=&
\sum_{i:t_i\leq t}  \int_{t_{i-1}}^{t_i} \bigg(
\bigg[
\Delta_n^{1/2}\sum_{k=i+1}^na'_{t_{k-1}}f(\Delta_n^{-1/2}\Delta^n_kW)1_{\{t_k\leq t\}}D_rX_{t_{k-1}}
\bigg] \int_{r}^1 \beta'_sD_rX_sds \bigg)dr
\eeas
for $t\in\Pi^n$. 
\end{en-text}
By the Clark-Ocone representation formula, we have 
\beas 
f'(\Delta_n^{-1/2}\Delta^n_iW)
&=& 
\Delta_n^{-1/2}
\int_{t_{i-1}}^{t_i}a_{n,i}(s)dW_{s}
\eeas
with 
\beas 
a_{n,i}(s) &=& \Delta_n^{1/2}
\E\bigg[D_{s}\bigg(f'(\Delta_n^{-1/2}\Delta^n_iW)\bigg)\>\big|\>{\cal F}_{s}\bigg],
\eeas
and moreover 
\beas 
a_{n,i}(s)
&=&
\E\big[f''(\Delta_n^{-1/2}\Delta^n_iW)\>|\>{\cal F}_{s}\big]1_{(t_{i-1},t_i]}(s)
\\&=&
g_{s}\big(\Delta_n^{-1/2}(W_{s}-W_{t_{i-1}})\big)\>1_{(t_{i-1},t_i]}(s)
\eeas
with 
\beas 
g_r(z) &=& \int f''\bigg(z+\sqrt{\frac{t_i-r}{\Delta_n}}\>x\bigg)\phi(x;0,1)dx
\eeas
for $r\in(t_{i-1},t_i]$. 
Then obviously,  
\beas 
\sup_{{s\in(t_{i-1},t_i] \atop i=1,...,n} \atop n\in\bbN}\|a_{n,i}(s)\|_{9,p}<\infty
\eeas
for every $p>1$. 
In the same way, we see that 
\beas 
f(\Delta_n^{-1/2}\Delta^n_iW)
&=& 
\Delta_n^{-1/2}\int_{t_{i-1}}^{t_i}\alpha_{n,i}(s)dW_{s}
\eeas
with some predictable processes $\alpha_{n,i}(s)$ satisfying 
\beas 
\sup_{{s\in(t_{i-1},t_i] \atop i=1,...,n} \atop n\in\bbN}\|\alpha_{n,i}(s)\|_{10,p}<\infty
\eeas
for every $p>1$. By Lemma 5 of \cite{Yoshida2012},
\beas&&
\bigg\| 
\sum_{i:t_i\leq t}  \int_{t_{i-1}}^{t_i}
 a_{t_{i-1}} f'(\Delta_n^{-1/2}\Delta^n_iW)
\Delta_n^{1/2}\sum_{k=i+1}^na'_{t_{k-1}}f(\Delta_n^{-1/2}\Delta^n_kW)1_{\{t_k\leq t\}}D_rX_{t_{k-1}}dr\bigg\|_{\mathbb L^{9}}
\\&=&
\bigg\| 
\Delta_n
\sum_{i:t_i\leq t} 
\bigg[ a_{t_{i-1}} 
\bigg(\Delta_n^{-1/2}\int_{t_{i-1}}^{t_i}a_{n,i}(s_1)dW_{s_1}\bigg)
\\&&
\times\bigg(
\Delta_n^{1/2}\sum_{k=i+1}^n
\bigg\{\int_{t_{i-1}}^{t_i}
a'_{t_{k-1}}1_{\{t_k\leq t\}} 
\Delta_n^{-1}D_rX_{t_{k-1}}
dr\bigg\}
\Delta_n^{-1/2}\int_{t_{k-1}}^{t_k}\alpha_{n,i}(s)dW_{s}
\bigg)\bigg]
\bigg\|_{\mathbb L^{9}}
\\&=&
O(\Delta_n^{1/2})
\eeas
for $t\in\Pi^n$. 
Hence 
\beas
\sup_{n\in\bbN}\sup_{t\in\Pi^n}\|\sigma_{11}(n,t)-\tilde{\sigma}_{11}(n,t)\|_{\mathbb L^{9}}
&=&
O(\Delta_n^{1/2})
\eeas
as $n\to\infty$. 
Furthermore, by the same lemma, we have 
\beas&&
\sup_{t\in\Pi^n}
\bigg\|
\sum_{i:t_i\leq t}  \int_{t_{i-1}}^{t_i} \bigg(
 a_{t_{i-1}} f'(\Delta_n^{-1/2}\Delta^n_iW)
 \int_{r}^1 \beta'_sD_rX_sds \bigg ) dr\bigg\|_{\mathbb L^{10}}
 \\&=&
 \sup_{t\in\Pi^n}
 \bigg\|
\Delta_n
\sum_{i:t_i\leq t}   a_{t_{i-1}} 
\bigg(\Delta_n^{-1/2}\int_{t_{i-1}}^{t_i}a_{n,i}(s_1)dW_{s_1}\bigg)
\bigg(\Delta_n^{-1}\int_{t_{i-1}}^{t_i} \bigg(\int_{r}^1 \beta'_sD_rX_sds \bigg) dr\bigg)\bigg\|_{\mathbb L^{9}}
\\&=&
O(\Delta_n^{1/2}).
\eeas
Therefore 
\beas
\sup_{t\in\Pi^n}\|\sigma_{12}(n,t)-\tilde{\sigma}_{12}(n,t)\|_{\mathbb L^{9}}
&=&
O(\Delta_n^{1/2}).
\eeas
{From} these estimates, 
\beas 
\sup_{t\in\Pi^n}\|\sigma(n,t)-\tilde{\sigma}(n,t)\|_{\mathbb L^{9}}
&=&
O(\Delta_n^{1/2}).
\eeas
One has
\bea\label{241228-1}
\det\tilde{\sigma}(n,t)
&=&
\tilde{\sigma}_{11}(n,t) \sigma_{22}(t)-\tilde{\sigma}_{12}(n,t)^2 
\nn\\&\geq&
\Delta_n\sum_{i:t_i\leq t}  \big[ a_{t_{i-1}} f'(\Delta_n^{-1/2}\Delta^n_iW)\big]^2
\times \sigma_{22}(t)
\nn\\&\geq&
\inf_x|a(x)|^2\>m_n\>\sigma_{22}(t)
\eea
for $t\in\Pi^n$, where the random variable $m_n$ is defined in Section 3.2.

Now, we shall verify $(C2^\flat)$. 
Checking $(C2^\flat)$ (iii) is not difficult if one estimates 
the $\mathbb H^{\otimes m}$-norms of $D_{r_1,....,r_m}$-derivative of the objects, 
in part with the aid of the Burkholder inequality.  

For $(C2^\flat)$ (ii), it suffices to show 
\bea\label{241227-3}
\limsup_{n\to\infty} \E[1_{\{m_n\geq {\sf c_1}\}}(\det \tilde{\sigma}(n,1/2))^{-p}] &<& \infty
\eea
for every $p>1$ since $s_n\geq 1/2$ when $m_n< {\sf c}_1$. 
Consider the two-dimensional stochastic process $\bar{X}_t=(X^{(1)}_t,X^{(2)}_t)$ 
defined by the stochastic integral equations with smooth coefficients 
\bea\label{241227-1}
\bar{X}_t  &=& \bar{X}_0+\int_0^t V_1(\bar{X}_s)\circ dW_s+\int_0^t V_0(\bar{X}_s)ds,
\eea
for $t\in[0,1]$, where the first integral is given in the Stratonovich sense and 
\beas 
V_1(x) =
\l[\begin{array}{c} b^{[1]}(x^1)\\0\end{array}\r],&&\sskip
V_0(x) = 
\l[\begin{array}{c} \tilde{b}^{[2]}(x^1)\\ \beta(x^1) \end{array}\r]
\eeas
for $x=(x^1,x^2)$,  
$\tilde{b}^{[2]}=b^{[2]}-2^{-1}b^{[1]}(b^{[1]})'$. 
\begin{en-text}
[Recall $\beta(x^1)=\mbox{Var}[f(Z)]a(x^1)^2=\mbox{Var}[f(Z)]|b^{[1]}(x^1)|^p$ and 
$b^{[1]}$ is uniformly nondegenerate. ]
\end{en-text}
Under $(C2)$, the system (\ref{241227-1}) satisfies the H\"ormander condition
\beas 
\mbox{Lie}[V_0;V_1](x^1,0)=\bbR^2\sskip(\forall x^1\in\mbox{supp}{\cal L}\{X_0\} ),
\eeas
the Lie algebra generated by $V_1$ and $V_0$, 
and as a result, for any $t\in(0,1]$ and $p>1$, there exists a constant $K_p$ such that 
\beas
\sup_{{\bf v}\in\bbR^2:|{\bf v}|=1}
\mathbb P\bigg[{\bf v}^\star \int_0^t \bar{Y}_s^{-1}V(\bar{X}_s)V(\bar{X}_s)^\star (\bar{Y}_s^{-1})^\star ds\>{\bf v}\leq\ep\bigg]
&\leq& K_p\ep^p
\eeas
for all $\ep\in(0,1)$. 
Here $\bar{Y}_t$ denotes a unique solution of the variational equation corresponding to (\ref{241227-1}). 
See Kusuoka and Stroock \cite{KusuokaStroock1984,KusuokaStroock1985}, 
Ikeda and Watanabe \cite{IkedaWatanabe1989}, Nualart \cite{N} 
for the nondegeneracy argument. 
Since both $\bar{Y}_1$ and $\bar{Y}_1^{-1}$ are bounded in $\cap_{p>1}\mathbb L^p$, we have 
\beas
\sup_{{\bf v}\in\bbR^2:|{\bf v}|=1}
\mathbb P\bigg[{\bf v}^\star \int_0^t \bar{Y}_1\bar{Y}_s^{-1}V(\bar{X}_s)V(\bar{X}_s)^\star (\bar{Y}_s^{-1})^\star \bar{Y}_1^\star ds\>{\bf v}\leq\ep\bigg]
&\leq& K_p'\ep^p
\eeas
form some constant $K_p'>0$, and in particular this implies 
\beas\label{241227-2}
\mathbb P\big[\sigma_{22}(t)\leq\ep\big]
&\leq& K_p'\ep^p
\eeas
for all $\ep\in(0,1]$. 
This inequality gives 
\beas 
\sigma_{22}(t)^{-1}&\in&\bigcap_{p>1} \mathbb L^p
\eeas
for every $t\in(0,1]$, and 
consequently, in view of (\ref{241228-1}), 
we obtained (\ref{241227-3}) 
and hence $(C2^\flat)$ (ii) for arbitrary ${\sf c}_1>0$. 

Finally, 
\beas &&
\sup_{t\geq\half} \mathbb P[\det\sigma_{(M^n_t,C_1)}<s_n]
\\&\leq&
\sup_{t\geq\half} \mathbb P[\det\sigma(n,t)<s_n]
\\&\leq&
\sup_{t\in\Pi^n:t\geq\half} \mathbb P[\det\sigma(n,t)<1.5s_n]
+\sup_{s,t:|t-s|\leq\Delta_n} \mathbb P\bigg[|\det\sigma(n,t)-\det\sigma(n,s)|>0.5s_n\bigg]
\\&\leq&
\mathbb P[\det\sigma(n,1/2)<1.5s_n]+O(\Delta_n^{1.35})
\\&\leq&
\mathbb P[\det\tilde{\sigma}(n,1/2)<2s_n]
+\mathbb P[|\det\sigma(n,1/2)-\det\tilde{\sigma}(n,1/2)|>0.5s_n]
+O(\Delta_n^{1.35})
\\&\leq&
\mathbb P[m_n>2{\sf c}_1,\>\det\tilde{\sigma}(n,1/2)<2s_n]+\mathbb P[m_n\leq 2{\sf c}_1]
\\&&
+\Delta_n^{-3/19}\E[|\det\sigma(n,1/2)-\det\tilde{\sigma}(n,1/2)|^3]
+2^{5\times19/3}\Delta_n^{5/3}\E[s_n^{-5\times19/3}]+O(\Delta_n^{1.35})
\\&=&
O(\Delta_n^{51/38})
\eeas
as $n\to\infty$ if we take ${\sf c}_1<\E[f'(Z)^2]/2$. 
Thus we have verified $(C2^\flat)$ (i), which completes the proof.  
\qed

\section{Stochastic expansion of generalized power variation of diffusions} \label{sec4} 
\setcounter{equation}{0}

\noindent
Hereafter we will concentrate on the stochastic expansion of the type \eqref{vndec} for the class
of generalized power variation. The results of this section are necessary 
for the derivation of the second order Edgeworth expansion for power variation, which is presented in Section \ref{241228-10},
but they might be also useful for other expansion problems in high frequency framework. 
We again consider a one-dimensional diffusion process $X=(X_t)_{t\in[0,1]}$ satisfying the
stochastic differential equation 
\begin{equation*} 
dX_t =  b^{[1]}(X_t) dW_t + b^{[2]}(X_t) dt.
\end{equation*}
Our aim is to study the stochastic expansion of generalized power variations of the form
\begin{equation} \label{vnf}
V_n(f) = \De_n \sum_{i=1}^{1/\De_n} f \Big( \frac{\De_i^n X}{\sqrt{\De_n}} \Big), \qquad \De_i^n X = X_{t_i} - X_{t_{i-1}},  
\end{equation}
where $f: \R \rightarrow \R$ is a given \textit{even} function, i.e. $f(x)=f(-x)$ for all $x \in \R$. This type
of functionals play a very important role in mathematical finance, where they are used for various estimation and testing
procedures; see e.g. \cite{BGJPS}, \cite{BS}, \cite{DPV} and \cite{J2} among many others. The most classical subclass of statistics (\ref{vnf}) are power variations, which correspond to functions
of the form $f(x)=|x|^p$; we will concentrate on Edgeworth expansion of power variations in the next section.
We introduce the notation
\begin{equation} \label{rhof}
\rho_x(f) = \mathbb E[f(xZ)], \qquad x\in \R, ~Z\sim \mathcal N(0,1) 
\end{equation}
whenever the latter is finite. Now, let us recall the law of large numbers and the central limit theorem for the functional
$V_n(f)$ derived in \cite{BGJPS}. 
\begin{en-text}
Below, we write $C^k_p (\R)$ to denote
the space of $k$-times continuously differentiable functions $h$ such that $h,\ldots, h^{(k)}$ have polynomial growth. 
\end{en-text}
\begin{theo} \label{th2}
(i) Assume that $b^{[1]},b^{[2]} \in C(\R )$ and $f\in C_p (\R)$. Then it holds that
\begin{equation} \label{lln}
V_n(f) \toop V(f)= \int_0^1 \rho_{b^{[1]}_s} (f) ds. 
\end{equation} 
(ii) If moreover $b^{[1]}\in C^{2}(\R )$ and $f\in C_p^1 (\R)$ we obtain the stable convergence 
\begin{equation} \label{clt}
\De_n^{-1/2} \Big(V_n(f) - V(f) \Big) \stab M \sim MN\left(0, \int_0^1 \rho_{b^{[1]}_s} (f^2) - \rho_{b^{[1]}_s}^2 (f) ds\right). 
\end{equation}
\end{theo}

\begin{rem} \rm
Recall that due to the It\^o formula the assumption $b^{[1]}\in C^{2}( \R )$ implies that the process 
$b^{[1]}_t$ satisfies a SDE of the form (\ref{sde}). Thus, $b^{[1]}_t$ is an It\^o semimartingale, which is usually
required for proving (\ref{clt}) (see e.g. \cite{BGJPS}). \qed  
\end{rem}
Now, we derive the second order stochastic expansion associated with the central limit theorem (\ref{clt}). 
Let us introduce the notation 
\begin{equation} \label{alpha}
\alpha_i^n = \De_n^{-1/2} b^{[1]}_{t_{i-1}} \Delta_i^n W,
\end{equation}
which serves as an approximation of the increment $\Delta_i^n X / \sqrt{\De_n}$.
One of the main results of this section  is the following theorem. We remark
that this result might be of independent  interest for other expansion problems in probability and statistics. 

\begin{theo} \label{th3}
Assume that $b^{[2]}\in C^{2}(\R )$, $b^{[1]}\in C^{4}(\R )$ and   
$f\in C_p^2 (\R)$. Then we obtain the stochastic expansion 
\begin{equation} \label{expansion}
 \widetilde{V}_n(f):=
\De_n^{-1/2} \Big(V_n(f) - V(f) \Big) = M_n + \De_n^{1/2} N_n + o_{\mathbb P}(\De_n^{1/2})
\end{equation}
with 
\begin{eqnarray} \label{Mn}
M_n= \De_n^{1/2} \sum_{i=1}^{1/\De_n} \Big(f ( \alpha_i^n) - \rho_{b^{[1]}_{t_{i-1}}} \Big),
\end{eqnarray} 
and $N_n=\sum_{k=1}^5 N_{n,k}$ 
\begin{eqnarray}
N_{n,1} &=& \De_n^{1/2} \sum_{i=1}^{1/\De_n} f'  (\alpha_i^n) \Big(b^{[2]}_{t_{i-1}} + \frac 12 b^{[1.1]}_{t_{i-1}}
H_2 (\Delta_i^n W / \sqrt{\Delta_n} ) \Big), \nonumber \\[1.5 ex]
N_{n,2}&=& \De_n^{-1/2} \sum_{i=1}^{1/\De_n}   f'(\alpha_i^n) \Big( b_{t_{i-1}}^{[2.1]} \int_{t_{i-1}}^{t_{i}} (W_s -W_{t_{i-1}}) ds \nonumber \\[1.5 ex]
&+& b_{t_{i-1}}^{[1.2]} \int_{t_{i-1}}^{t_i}  \{s-t_{i-1} \} dW_s  
+\frac{\Delta_n^{3/2}b_{t_{i-1}}^{[1.1.1]}}{6}   H_3 (\Delta_i^n W / \sqrt{\Delta_n} )  \Big),     \nonumber \\[1.5 ex]
N_{n,3}&=& \frac{\De_n}{2} \sum_{i=1}^{1/\De_n}   f''(\alpha_i^n) \Big( b_{t_{i-1}}^{[2]} + \frac 12 b_{t_{i-1}}^{[1.1]}
H_2 (\Delta_i^n W / \sqrt{\Delta_n} ) \Big)^2, \nonumber \\[1.5 ex]
N_{n,4}&=& \frac{1}{2\De_n} \sum_{i=1}^{1/\De_n} \Big(  -\rho''_{b_{t_{i-1}}^{[1]}} (f) 
|b_{t_{i-1}}^{[1.1]}|^2 \int_{t_{i-1}}^{t_{i}} (W_s -W_{t_{i-1}})^2 ds  \nonumber \\[1.5 ex]
&-& \Delta_n^2 \rho'_{b_{t_{i-1}}^{[1]}} (f) b_{t_{i-1}}^{[1.2]} \Big), \nonumber \\[1.5 ex]
N_{n,5}&=&  -\Delta_n^{-1}\sum_{i=1}^{1/\De_n} \rho'_{b_{t_{i-1}}^{[1]}} (f) b_{t_{i-1}}^{[1.1]} 
\int_{t_{i-1}}^{t_{i}} (W_s -W_{t_{i-1}})ds,
\end{eqnarray}
where $(H_k)_{k\geq 0}$ denote the Hermite polynomials 
and the processes $b_t^{[k_1 \ldots k_d]}$ were defined in Section \ref{sec3}.
\end{theo}
{\it Proof.} See Section \ref{proof}. \qed \newline \newline
To describe the limits of the quantities $N_{n,k}$, $1\leq k\leq 5$, we need to introduce some further notation.
\newline \newline
{\it Notation.} We introduce the functions $g_k: \R^6 \rightarrow \R$, $1\leq k\leq 5$, as follows:
\begin{eqnarray*}
g_1(x_1,\ldots, x_6)&=& \E \Big[U f'  (x_2 U) \Big(x_1 + \frac 12 x_5
H_2 (U) \Big) - \rho'_{x_2}(f) x_5 UV \Big]\\[1.5 ex] 
g_2(x_1,\ldots, x_6)&=& \E \Big[ f'(x_2 U)\Big((x_3+x_4) V +  \frac 16 x_6 H_3(U)\Big) \Big] \\[1.5 ex] 
g_3(x_1,\ldots, x_6)&=&   \frac 12 \E \Big[ f''(x_2 U) \Big( x_1 + \frac 12 x_5
H_2 (U) \Big)^2 \Big] \\[1.5 ex] 
g_4(x_1,\ldots, x_6)&=&  - \frac{1}{4} \rho''_{x_2} (f) 
x^2_5 - \frac 12 \rho'_{x_2} (f) x_4 \\[1.5 ex]
g_5(x_1,\ldots, x_6)&=& \E \Big[ \Big \{ f'  (x_2 U) \Big(x_1 + \frac 12 x_5
H_2 (U) \Big) - \rho'_{x_2}(f) x_5 V \Big \}^2\Big]  
\end{eqnarray*}
with 
$$(U,V)\sim \mathcal N_2 \left(0, 
\Big( \begin{array} {cc}
1 & 1/2 \\
1/2 & 1/3
\end{array}
\Big)
\right).$$

\begin{rem} \rm
Theorem \ref{generalucp} implies the convergence in probability
\begin{equation} \label{nklim}
N_{n,k} \toop N_{k}=\int_0^1 g_k(b^{[2]}_s, b^{[1]}_s, b_s^{[2.1]}, b^{[1.2]}_s, b^{[1.1]}_s,b^{[1.1.1]}_s) ds, \qquad k=2,3,4
\end{equation}
under the assumptions of Theorem \ref{th3}. The terms $N_{n,1}$ and $N_{n,5}$ converge stably in law due to Theorem 
\ref{generalclt}; their asymptotic distributions will be specified later. \qed    
\end{rem}

\begin{rem} \rm
The fact that we consider the drift and volatility processes of the type $b_s^{[k]} = b^{[k]}(X_s)$ is not essential
for developing the stochastic expansion of Theorem \ref{th3}. In general the processes $b_s^{[k_1 \ldots k_l]}$ 
that appear in Theorem \ref{th3} may depend on different Brownian motions, which are not perfectly correlated with 
$W$ that drives the process $X$. In this case a similar stochastic expansion can be deduced; however, it will contain 
additional terms, which are due to new Brownian motions. \qed    
\end{rem}
In the next section we will require a consistent estimator of the asymptotic variance of $M_n$, i.e.
\begin{equation*} 
C= \int_0^1 \rho_{b^{[1]}_s} (f^2) - \rho_{b^{[1]}_s}^2 (f) ds.
\end{equation*}
A rather natural one is given by
\begin{equation} \label{fnc} 
F_n = \Delta_n \sum_{i=1}^{1/\De_n} f^2 \Big( \frac{\De_i^n X}{\sqrt{\De_n}} \Big) -
f \Big( \frac{\De_i^n X}{\sqrt{\De_n}} \Big) f \Big( \frac{\De_{i+1}^n X}{\sqrt{\De_n}} \Big)   
\end{equation} 
We remark that $F_n$ is a feasible statistic in contrast to the Riemann sum approximation defined at \eqref{riemannfn}. 
The next theorem, which follows from the combination of central limit theorems
presented in \cite{BGJPS} and Theorem \ref{generalclt}, describes the joint asymptotic distribution of $(M_n, F_n, N_n)$.
This result is crucial for the derivation of the second order Edgeworth expansion.

\begin{theo} \label{th5}
Assume that conditions of Theorem \ref{th3} are satisfied. Then we obtain the stable convergence
\begin{equation*} 
\Big(M_n,  \Delta_n^{-1/2}(F_n -C), N_n \Big) \stab 
(M, \widehat{F}, N) \sim MN\Big(\mu, \int_0^1 \Xi_s ds \Big),  
\end{equation*}
where the matrix $\Xi_s$ is given as
\begin{eqnarray*}
\Xi_s^{11} &=& \rho_{b^{[1]}_s} (f^2) - \rho_{b^{[1]}_s}^2 (f) \\[1.5 ex]
\Xi_s^{12} &=& \Xi_s^{21} = \rho_{b^{[1]}_s} (f^3) - 3\rho_{b^{[1]}_s} (f^2) \rho_{b^{[1]}_s} (f) 
+2\rho_{b^{[1]}_s}^3 (f) \\[1.5 ex]
\Xi_s^{22} &=& \rho_{b^{[1]}_s} (f^4) - 4 \rho_{b^{[1]}_s} (f^3) \rho_{b^{[1]}_s} (f) 
+6 \rho_{b^{[1]}_s} (f^2) \rho_{b^{[1]}_s}^2 (f) -3\rho_{b^{[1]}_s}^4 (f), \\[1.5 ex]
\Xi_s^{33} &=& (g_5 -g_1^2)(b^{[2]}_s, b^{[1]}_s, b_s^{[2.1]}, b^{[1.2]}_s, b^{[1.1]}_s,b^{[1.1.1]}_s), 
\end{eqnarray*}
and $\Xi_s^{13} = \Xi_s^{23}=0$, and $\mu_1=\mu_2=0$,
\begin{equation*} 
\mu_3= \int_0^1 g_1(
 b^{[2]}_s, b^{[1]}_s, b_s^{[2.1]}, b^{[1.2]}_s, b^{[1.1]}_s,b^{[1.1.1]}_s
) dW_s + \sum_{k=2}^4 N_k.  
\end{equation*}
\end{theo}

\section{Asymptotic expansion for the power variation}
\label{241228-10}
\setcounter{equation}{0}
%
Now we have all instruments at hand to obtain the Edgeworth expansion for the case
of power variation $V_n(f_p)$ with
\begin{equation*}
f_p(x) = |x|^{p},   
\end{equation*}
which is our leading example. 
As we mentioned in Section \ref{sec4}, this would be 
the most important class of functionals in mathematical finance.
In order to obtain the Edgeworth expansion for power variation, we will combine the results of Sections \ref{sec3} and \ref{sec4}.
Applying Theorem \ref{th3} to the function $f_p$ we see that the martingale part $M_n$ is given as 
\begin{eqnarray*}
M_n= \De_n^{1/2} \sum_{i=1}^{1/\De_n} |b^{[1]} (X_{t_{i-1}})|^{p} 
\Big( \Big| \frac{\De_i^n W}{\sqrt{\De_n}} \Big|^{p} - m_{p} \Big)
\end{eqnarray*}
with $m_{p} = \E[|\mathcal N(0,1)|^{p}]$. In particular, $M_n$ is a weighted power variation studied in Section \ref{sec3}.
Consequently, we can apply the results of Section \ref{sec3} with 
\begin{equation*}
a(x)= |b^{[1]} (x)|^{p}, \qquad f(x)= f_{p} (x) -m_{p} \qquad \mbox{and} \quad p\in 2\N \cup  (11,\infty). 
\end{equation*}
Now, we will compute all quantities from previous sections required for the Edgeworth expansion.   
First, we obtain the Hermite expansion 
\begin{equation*}
f(x) = \sum_{k=2}^{\infty} \lambda_k H_k(x)
\end{equation*} 
with $\lambda_k=0$ if $k$ is odd (because $f$ is an even function), and
\begin{equation*}
\lambda_2 =\frac{m_{p+2} - m_{p}}{2}.
\end{equation*}
We start with the computation of the random symbol $\underline{\sigma}$. Here we mainly need to determine the functions 
$g_1,\ldots,g_5$ defined in Section \ref{sec4}. We observe that, for any $k\geq 0$ with $k<p$,
\begin{equation*}
f^{(k)}_p(x) = \mbox{sgn}(x)^k p(p-1)\cdots (p-k+1) |x|^{p-k}, \qquad \rho_x(f_p)= m_p |x|^p. 
\end{equation*} 
Now, a straightforward calculation gives the identities
\begin{eqnarray*}
g_1(x_1,\ldots, x_6)&=& p~ \mbox{sgn}(x_2) |x_2|^{p-1} \Big(x_1 m_p + \frac 12 x_5 
 (m_{p+2}-2m_p) 
\Big) \\[1.5 ex] 
g_2(x_1,\ldots, x_6)&=& p~ \mbox{sgn}(x_2) |x_2|^{p-1} \Big(\frac 12 (x_3+x_4)m_p + \frac 16 x_6(m_{p+2} -m_p) \Big) \\[1.5 ex] 
g_3(x_1,\ldots, x_6)&=&   \frac {p(p-1)}{2} |x_2|^{p-2} 
\Big(x_1^2 m_{p-2} + x_1x_5(m_{p} -m_{p-2}) + \frac {x_5^2}{4}(m_{p+2} - 2m_p + m_{p-2}) \Big)\\[1.5 ex] 
g_4(x_1,\ldots, x_6)&=&  \frac {p}{4} m_p 
\Big( -(p-1)|x_2|^{p-2} x_5^2 -2 x_4~\mbox{sgn}(x_2) |x_2|^{p-1} \Big) \\[1.5 ex]
g_5(x_1,\ldots, x_6)&=& p^2 |x_2|^{2p-2} \Big(x_1^2 m_{2p-2} + x_1x_5 (m_{2p} - m_{2p-2}) + 
\frac{x_5^2}{4} (m_{2p+2} - 2m_{2p} + m_{2p-2}) \\[1.5 ex]
 &+& \frac{x_5^2}{3} m_p^2 - x_5 m_p (x_1 m_p + \frac{x_5}{2} [m_{p+2} - m_p]) \Big)
\end{eqnarray*}
As in the previous section we consider the quantity
\begin{equation*} 
F_n = \Delta_n \sum_{i=1}^{1/\De_n} f_{2p} \Big( \frac{\De_i^n X}{\sqrt{\De_n}} \Big) -
f_{p} \Big( \frac{\De_i^n X}{\sqrt{\De_n}} \Big) f_{p} \Big( \frac{\De_{i+1}^n X}{\sqrt{\De_n}} \Big)   
\end{equation*} 
as a consistent estimator of $C$. We obtain the following result, which again follows from Theorem \ref{generalclt}. 

\begin{theo} \label{th6}
Assume that conditions of Theorem \ref{th3} are satisfied. Then we obtain the stable convergence
\begin{equation*} 
\Big(M_n,  \Delta_n^{-1/2}(F_n -C), N_n, \Delta_n^{-1/2}(C_n -C) \Big) \stab 
(M, \widehat{F}, N,\widehat{C})\sim MN\Big(\mu, \int_0^1 \Xi_s ds \Big),  
\end{equation*}
where the entries $\Xi_s^{ij}$, $1\leq i,j\leq 3$, of the matrix $\Xi_s\in \R^{4 \times 4}$ and $\mu_j$, $1\leq j\leq 3$ 
of the vector $\mu \in \R^4$ are given in Theorem \ref{th5},
and $\mu_4=\Xi_s^{34}=0$,
\begin{eqnarray*}  
\Xi_s^{14}&=&\Xi_s^{41}= \Gamma_2~ |b^{[1]} (X_s)|^{3p}, \\[1.5 ex]
\Xi_s^{24}&=&\Xi_s^{42}= \overline{\Gamma}~ |b^{[1]} (X_s)|^{4p}, \\[1.5 ex]
\Xi_s^{44}&=& \Gamma_1~ |b^{[1]} (X_s)|^{4p}, 
\end{eqnarray*} 
where the constants $\Gamma_1, \Gamma_2$   are given in Proposition  \ref{prop3} and  $\overline{\Gamma} $ is defined as
\begin{equation*} 
\overline{\Gamma} = \mbox{Cov} 
\left[ f_{2p}( W_1 ),
\int_{0}^{1} \E^2[f'_p(W_1)| \mathcal F_s] ds \right] - 2
 \mbox{Cov} 
\left[ f_p( W_1 ) f_p( W_2-W_1 ),
\int_{0}^{1} \E^2[f'_p(W_1)| \mathcal F_s] ds \right].
\end{equation*}
\end{theo}
As a consequence of Theorem \ref{th6} and Remark \ref{rem1} we conclude that  
\bea\label{241231-21}
\underline \sigma (z, \iu,\iv) 
&=& 
(\iu)^2 \mathcal{H}_1(z)  +  \iu \mathcal{H}_2 + \iv\mathcal{H}_3(z)
\eea
with 
\begin{equation*}
\mathcal{H}_1(z) = z\frac{\int_0^1 \Xi_s^{14} ds}{2\int_0^1 \Xi_s^{11} ds}, \qquad \mathcal{H}_2 = \mu_3, \qquad 
\mathcal{H}_3(z) = z \frac{\int_0^1 \Xi_s^{12} ds}{\int_0^1 \Xi_s^{11} ds}.  
\end{equation*}
It should be noted that $\underline{\sigma}$ of (\ref{241231-21}) 
is essentially the same but different from $\underline{\sigma}$ of (\ref{241231-2}) since 
the reference functional $F_n$ is now defined by (\ref{fnc}) not by (\ref{241231-10}) 
while the limits of both coincide each other and the ways of derivation of two adaptive random symbols 
are the same except for $\widehat{F}$.  
Using the results of Section \ref{sec3} we immediately obtain the anticipative random symbol
\bea\label{241231-22}
\overline \sigma (\iu,\iv) = \iu \Big(\iv- \frac{u^2}{2} \Big)^2 \mathcal{H}_4 + \iu \Big(\iv- \frac{u^2}{2} \Big) \mathcal{H}_5 
\eea
with 
\begin{equation*} 
\mathcal{H}_4 = \lambda_2 (m_{2p} - m_p^2)^2 \mathcal C_2, \qquad 
\mathcal{H}_5 =\lambda_2 (m_{2p} - m_p^2)(\mathcal C_3 + \mathcal C_4),
\end{equation*} 
where 
\begin{eqnarray*} 
\mathcal  C_2 &=&  \int_0^1 |b^{[1]}(X_{s})|^{p}
  \Big(\int_{s}^1 \Big(|b^{[1]}|^{2p} \Big)'(X_u) D_{s} X_u du \Big)^2 ds, \\[1.5 ex]
\mathcal  C_3 &=&  
\int_0^1 |b^{[1]}(X_{s})|^{p} \Big( \int_{s}^1 \Big(|b^{[1]}|^{2p} \Big)'' (X_u) (D_s X_u)^2 du  \Big) ds, \\[1.5 ex]
\mathcal  C_4 &=&   \int_0^1 |b^{[1]} (X_{s})|^{p} \Big( \int_{s}^1 \Big(|b^{[1]}|^{2p} \Big)'(X_u) D_s D_s X_u du  \Big) ds.
\end{eqnarray*}
In the power variation case, $a(x)=|b^{[1]}(x)|^p$ and we assumed in $(C1)$ that 
$a(x)$ is bounded away from zero. 
So, in our situation, $a(x)$ is smooth in a neighborhood of $X_0$. 
By a certain large deviation argument, 
we may assume that $a(x)$ is smooth and even having bounded derivatives, 
from the beginning, at least in the proof of asymptotic nondegeneracy. 

{From} the above argument, we obtain an asymptotic expansion for the power variation. 
Recall $\widetilde{V}_n(f)=\De_n^{-1/2} \big(V_n(f) - V(f) \big)$. 

\begin{theo} \label{th52}
Let $b^{[1]},b^{[2]}\in C^\infty_{b,1}(\bbR)$ and $f_p(x)=|x|^p$ with 
$p\in 2\N \cup  ( 13,\infty)$.
Assume that $\inf_{x} |b^{[1]}(x)|>0$, $\sum_{k=1}^\infty|(b^{[1]})^{(k)}(X_0)|>0$ and  
let the functional $F_n$ be given by (\ref{241231-10}). 
Then for the density $p_n(z,x)$ corresponding to the random symbol $\sigma$ 
determined by (\ref{241231-21}) and (\ref{241231-22}), it holds that 
\beas 
\sup_{h\in\mathcal{E}(K,\gamma)}
\bigg|\E[f(\widetilde{V}_n(f_p),F_n)]-\int h(z,x)p_n(z,x)dzdx\bigg|
&=&
o\big(\sqrt{\Delta_n}\>\big)
\eeas
as $n\to\infty$, for any positive numbers $K$ and $\gamma$.
\end{theo}
Theorem \ref{th52} is proved by applying Theorems \ref{241227-8} and \ref{th6}. 
In the present situation, $N_n$ involves $f''$ and that is the reason why the number 13 appears. 
However, it would be possible to reduce it to 11 if the estimations related with $N_n$-part is refined, 
though we do not pursue this point in this article.

Theorem \ref{th52} and the corresponding Edgeworth expansion for the studentized statistics at \eqref{edex} are the main 
results of this paper. In particular, these asymptotic expansions can be applied to distribution analysis of various statistics
in financial mathematics as power variation type estimators are frequently used in this field. Another potential area of application
is Euler approximation of continuous SDE's of the form \eqref{sde}. As is well-known from \cite{JP}, the Euler approximation scheme 
is asymptotically mixed normal and its limit depends on the asymptotic theory for quadratic variation. Thus, our Edgeworth expansion
results can be potentially applied to numerical analysis of SDE's to obtain a more precise formula for the error distribution.

\begin{rem} \rm
As we mentioned above, we can combine the results of Sections \ref{sec3} and \ref{sec4}, because we consider the power
function $f_p(x)=|x|^p$. In this case the dominating part $M_n$ is a weighted power variation in the sense of 
Section \ref{sec3}.
The case of a general {\it even} function $f$ is more complicated. The results of \ref{sec4} still apply, but the computations
of the random symbol $\overline \sigma$ is more involved. Let us shortly sketch the idea how $\overline \sigma$ can be obtained. 
Recall that in the general case the term $M_n$ is given as 
\begin{eqnarray*} 
M_n= \De_n^{1/2} \sum_{i=1}^{1/\De_n} \Big(f ( \alpha_i^n) - \rho_{b^{[1]}_{t_{i-1}}} \Big)
\end{eqnarray*} 
(see Theorem \ref{th3}). As in Section \ref{sec3} we therefore need to compute the projection onto the
second order Wiener chaos of the quantity
\begin{equation*}
f \Big( \De_n^{-1/2} b^{[1]}_{t_{i-1}} \Delta_i^n W \Big) - \rho_{b^{[1]}_{t_{i-1}}}.
\end{equation*} 
For this purpose we use the following multiplication formula (see \cite{AS})
\begin{equation*}
H_k(\gamma x) =\sum_{i=0}^{[k/2]} 2^{-i}\gamma^{k-2i} (\gamma ^2 - 1)^{i} {k \choose 2i} 
\frac{(2i)!}{i!} H_{k-2i} (x), \qquad \gamma \in \R.
\end{equation*}
Under the assumptions of Section \ref{sec4}, the function $f$ admits the Hermite expansion 
$f(x)= \sum_{k=0}^{\infty} \lambda_{2k} H_{2k} (x)$ (since $f$ is even). Hence, we deduce that
\begin{equation*}
f(\gamma x) = \sum_{k=0}^{\infty} \lambda_{2k} \left ( 
\sum_{i=0}^{k} 2^{-i}\gamma^{2k-2i} (\gamma ^2 - 1)^{i} {2k \choose 2i} \frac{(2i)!}{i!} H_{2k-2i} (x) \right).
\end{equation*}
We conclude that the projection of $f ( \alpha_i^n) - \rho_{b^{[1]}_{t_{i-1}}}$ onto the second order Wiener chaos
is given by 
\begin{equation*}
\sum_{k=1}^{\infty} 2^{-k+1} \lambda_{2k}  
|b^{[1]}_{t_{i-1}}|^{2} (|b^{[1]}_{t_{i-1}}|^2 - 1)^{k-1} {2k \choose 2} \frac{(2(k-1))!}{(k-1)!} 
~H_2 (\De_n^{-1/2}  \Delta_i^n W ).
\end{equation*} 
Using this identity one can compute $\overline \sigma$ as in Section \ref{sec3}. However, we dispense with the exact 
exposition. \qed  
\end{rem}

\section{Studentization} \label{sec6}
\setcounter{equation}{0}
As we mentioned in the beginning, we are mainly interested in the second order Edgeworth expansion connected with standard central limit theorem
\begin{equation*}
\frac{Z_n}{\sqrt{F_n}} \schw \mathcal N(0,1).
\end{equation*}
where $F_n$ is a consistent estimator of $C$ defined in (\ref{fnc}). In the following we present such an Edgeworth expansion 
for the case of power variation discussed in the previous section. First of all, we remark that the random symbol $\sigma(z,\iu,\iv)$
is given as
\begin{equation*}
\sigma(z,\iu,\iv) = \sum_{j=1}^8 c_j(z)(\iu)^{m_j} (\iv)^{n_j}, 
\end{equation*}
where
\begin{eqnarray*}
&&m_1=1, \quad n_1=0, \quad c_1(z)=\mathcal H_2, \qquad m_2=0, \quad n_2=1, \quad c_2(z)=\mathcal H_3(z) \\
&&m_3=2, \quad n_3=0, \quad c_3(z)=\mathcal H_1(z), \qquad m_4=1, \quad n_4=1, \quad c_4(z)=\mathcal H_5 \\ 
&&m_5=3, \quad n_5=0, \quad c_5(z)=\frac 12 \mathcal H_5, \qquad m_6=1, \quad n_6=2, \quad c_6(z)=\mathcal H_4 \\
&&m_7=3, \quad n_7=1, \quad c_7(z)=\mathcal H_4, \qquad m_8=5, \quad n_8=0, \quad c_8(z)=\frac 14 \mathcal  H_4. 
\end{eqnarray*}
As a consequence, we obtain the following decomposition for the density $p_n(z,x)$ of $(Z_n, F_n)$:
\begin{equation*}
p_n(z,x) =  \phi(z;0,x) p^C(x) + \Delta_n^{1/2} \sum_{j=1}^8 p_j(z,x) 
\end{equation*}
with
\begin{equation*}
p_j(z,x) =  (-d_z)^{m_j} (-d_x)^{n_j}
\Big( \phi(z;0,x)  p^C(x) \E \left[  c_j(z)|C=x\right]  \Big), \qquad 1\leq j\leq 8. 
\end{equation*}
We start with the following observation. Let $\Pi$ be a finite measure on $\R$ with density $\pi$, such that all moments
of $\Pi$ are finite. Then it trivially holds that
\begin{equation*}
\lim_{x \rightarrow \infty } |x|^k \pi(x) =0, \qquad \lim_{x \rightarrow -\infty } |x|^k \pi(x) =0 \qquad k\geq 0.   
\end{equation*} 
Given that the density $\pi$ is a $C^k$ function and $g$ is a polynomial, we also have
\begin{equation*}
\int_{\R} g^{(k)}(x) \pi(x) dx = (-1)^k \int_{\R} g(x) \pi^{(k)}(x) dx    
\end{equation*}
by induction. Let $g$ be an arbitrary polynomial and $\kappa (x)= \E[H|C=x] p^C(x)$ for an integrable random variable $H$, and note
that 
\begin{equation*}
\int_\R m(x) \kappa (x)dx = \E[m(C)H],
\end{equation*}
whenever the integral makes sense. 
We define the polynomials $q_{\beta,v} (z, x)$ via
\begin{equation*}
d_x^{\beta} g(z/\sqrt{x}) = \sum_{v\leq \beta} q_{\beta,v} (z/\sqrt{x}, 1/\sqrt{x}) g^{(v)} (z/\sqrt{x}),   
\end{equation*}
where $g^{(v)}$ denotes the $v$th derivative of $g$. Let $(\alpha,\beta )\in \N_0^2$. Then it holds that 

\begin{eqnarray*}
&&\int_{\R^2} g \Big( \frac{z}{\sqrt{x}}\Big) d_z^\alpha d_x^\beta \Big[\phi(z;0,x) \kappa (x) \Big] dz dx = 
(-1)^\beta \int_{\R^2} d_x^\beta g \Big( \frac{z}{\sqrt{x}}\Big)  d_z^\alpha  \phi(z;0,x)   \kappa (x) dz dx \\[1.5 ex]
&&= (-1)^\beta \int_{\R^2} \sum_{v\leq \beta} q_{\beta,v} \Big( \frac{z}{\sqrt{x}}, \frac{1}{\sqrt{x}}\Big) 
g^{(v)} \Big(\frac{z}{\sqrt{x}} \Big)  
d_z^\alpha  \phi(z;0,x)   \kappa (x) dz dx \\[1.5 ex]
&&= (-1)^\beta \int_{\R^2} \sum_{v\leq \beta} q_{\beta,v} \Big( y, \frac{1}{\sqrt{x}}\Big) 
g^{(v)} (y )  
x^{-\alpha /2}d_y^\alpha  \phi(y;0,1)   \kappa (x) dy dx \\[1.5 ex]
&&= (-1)^\beta \int_{\R} g(y) \sum_{v\leq \beta} (-1)^v d_y^v \left\{ d_y^\alpha  \phi(y;0,1) 
\int_{\R} q_{\beta,v} \Big( y, \frac{1}{\sqrt{x}}\Big)   
x^{-\alpha /2}   \kappa (x)  dx\right \} dy  \\[1.5 ex]
&& = \int_{\R} g(y) \sum_{v\leq \beta} (-1)^{\beta +v} d_y^v \left\{ d_y^\alpha  \phi(y;0,1) 
\E\Big[H C^{-\alpha /2} q_{\beta,v} ( y, C^{-1/2} ) \Big] \right \} dy.  
\end{eqnarray*}
Clearly, the above identity will enable us to compute the Edgeworth expansion for the studentized statistic $Z_n/\sqrt{F_n}$. We need
to determine the polynomials $q_{\beta,v}$ for $\beta=0,1,2$:
\begin{eqnarray*}
&& q_{0,0}(a,b)=1, \\
&& q_{1,0} (a,b)=0, \qquad  q_{1,1} (a,b)= -\frac 12 ab^2, \\
&& q_{2,0} (a,b)=0, \qquad  q_{2,1} (a,b)= \frac 34 ab^4, \qquad  q_{2,2} (a,b)= \frac 14 a^2 b^4.
\end{eqnarray*}
Recall the identity $d_y^\alpha  \phi(y;0,1) = (-1)^\alpha H_{\alpha} (y) \phi(y;0,1)$ and
\begin{equation*}
H_1(x)=x, \qquad H_3(x)=x^3-3x, \qquad H_5(x)=x^5-10x^3 + 15x.
\end{equation*}
A straightforward computation shows that 
\begin{eqnarray*}
\int_{\R^2} g \Big( \frac{z}{\sqrt{x}}\Big) p_1(z,x) dzdx &=& \E[\mathcal H_2 C^{-1/2}] \int_{\R} g(y) y \phi(y;0,1)  dy, \\[1.5 ex]
\int_{\R^2} g \Big( \frac{z}{\sqrt{x}}\Big) \sum_{j=4}^5 p_j(z,x) dzdx 
&=& -\frac 12 \E[\mathcal H_5 C^{-3/2}] \int_{\R} g(y) y \phi(y;0,1)  dy, \\[1.5 ex] 
\int_{\R^2} g \Big( \frac{z}{\sqrt{x}}\Big) \sum_{j=6}^8 p_j(z,x) dzdx 
&=& \frac{3}{4}
\E[\mathcal H_4 C^{-5/2}] \int_{\R} g(y) y \phi(y;0,1)  dy. 
\end{eqnarray*}  
The corresponding computation for the terms $p_2(z,x)$ and $p_3(z,x)$ has to be performed separately, since the random 
variables $c_2$ and $c_3$ depend on $z$. Recall that the quantities $\mathcal H_1(z)$ and $\mathcal H_3(z)$ are linear in $z$,
i.e. $\mathcal H_1(z)=z \widetilde{\mathcal H}_1$, $\mathcal H_3(z)=z \widetilde{\mathcal H}_3$. We deduce as above (here 
$\kappa (x)= \E[\widetilde{\mathcal H}_3|C=x] p^C(x)$)  
\begin{eqnarray*}
&& \int_{\R^2} g \Big( \frac{z}{\sqrt{x}}\Big) p_2(z,x) dzdx = -
\int_{\R^2} z g \Big( \frac{z}{\sqrt{x}}\Big)  d_x \Big[\phi(z;0,x) \kappa (x) \Big] dz dx \\[1.5 ex] 
&& = \int_{\R} g(y)  d_y \left\{ y \phi(y;0,1) 
\E\Big[\widetilde{\mathcal H}_3  q_{1,1} ( y, C^{-1/2} ) C^{1/2} \Big] \right \} dy \\[1.5 ex]
&& = \frac 12 \E[\widetilde{\mathcal H}_3  C^{-1/2}  ]
 \int_{\R} g(y)   \phi(y;0,1) (2y-y^3)  dy.
\end{eqnarray*} 
Finally, we obtain that (here 
$\kappa (x)= \E[\widetilde{\mathcal H}_1|C=x] p^C(x)$)  
\begin{eqnarray*}
&& \int_{\R^2} g \Big( \frac{z}{\sqrt{x}}\Big) p_3(z,x) dzdx = 
\int_{\R^2}  g \Big( \frac{z}{\sqrt{x}}\Big)  d_z^2 \Big[z \phi(z;0,x) \kappa (x) \Big] dz dx \\[1.5 ex] 
&& = \int_{\R^2} x^{-1} g''(y)   y \phi(y;0,1) \kappa (x) dy dx  \\[1.5 ex]
&& =  \E[\widetilde{\mathcal H}_1  C^{-1/2}  ]
 \int_{\R} g(y)   d_y^2 [y\phi(y;0,1)]  dy = \E[\widetilde{\mathcal H}_1  C^{-1/2}  ]
 \int_{\R} g(y)   H_3(y) \phi(y;0,1)   dy.
\end{eqnarray*}
Combining the above results, we deduce the second order Edgeworth expansion for the density of $Z_n/\sqrt{F_n}$
\begin{eqnarray}
&& \label{edex} p^{Z_n/\sqrt{F_n}} (y) = \phi(y;0,1) + \Delta_n^{1/2} \phi(y;0,1) \Big(y \Big\{ \E[\mathcal H_2 C^{-1/2}] 
 -\frac 12 \E[\mathcal H_5 C^{-3/2}]  \\[1.5 ex]
&&+  \frac{3}{4}
\E[\mathcal H_4 C^{-5/2}] + \E[\widetilde{\mathcal H}_3  C^{-1/2}  ] - 3 \E[\widetilde{\mathcal H}_1  C^{-1/2}  ]
 \Big\} + y^3 \Big\{ \E[\widetilde{\mathcal H}_1  C^{-1/2}  ] - \frac 12 \E[\widetilde{\mathcal H}_3  C^{-1/2}  ]\Big\}\Big),
\nonumber 
\end{eqnarray}
which is one of the main statements of the paper.

\begin{rem} \rm
In practice the application of the asymptotic expansion at \eqref{edex} requires the knowledge of the coefficients of the type
$b^{[k_1\ldots k_d]}$ (cf. \eqref{nklim}). While the volatility related processes $b^{[1]}$, $b^{[1.1]}$, $b^{[1.1.1]}$
can be estimated from high frequency data $X_{t_i}$, the drift related processes $b^{[2]}$, $b^{[2.1]}$, $b^{[1.2]}$
can't be consistently estimated on a fixed time span. Thus, the applicability of the Edgeworth expansion
at \eqref{edex} relies on the knowledge of the drift related coefficients or their estimation on an infinite time span. \qed   
\end{rem}

\begin{ex} \rm ({\it Classical Edgeworth expansion})
In this example we compare the classical Edgeworth expansion with the result derived in (\ref{edex}). Let $(Y_i)_{i\geq 1}$
be a sequence of i.i.d. random variables with mean $\mu$ and variance $\sigma^2$. Define $S_n=n^{-1/2} \sum_{i=1}^n
\sigma^{-1} (Y_i - \mu)$. Then the second order Edgeworth expansion of the density of $S_n$ is given as
\begin{equation*}
\phi(y;0,1) + \frac{\kappa_3}{6\sigma ^3 \sqrt{n}} \phi(y;0,1) H_3(y), 
\end{equation*}
where $\kappa_3$ denotes the third cumulant of the law of $Y_1$. Let us now consider the quantity $M_n$ from
(\ref{weightedw}) with $a\equiv 1$ and $\Delta_n=n^{-1}$, i.e. 
\begin{equation*}
M_n= n^{-1/2} \sum_{i=1}^{n} f(\sqrt n \Delta_i^n W) \qquad \mbox{with} \qquad \E[f(W_1)]=0. 
\end{equation*}  
Due to self-similarity of the Brownian motion, we are in the classical setting of i.i.d. observations.
In this case $C=\sigma^2=\E[f^2(W_1)]$, and if we set $F_n\equiv C$, we obtain from (\ref{edex}): 
\begin{equation*}
\phi(y;0,1) +  \frac{1}{\sqrt n}\phi(y;0,1) H_3(y) \widetilde{\mathcal H}_1 C^{-1/2}  
\end{equation*} 
as the approximative density, since all quantities in (\ref{edex}) are $0$ except $C$ and $\widetilde{\mathcal H}_1$. 
We now show that the quantities $\widetilde{\mathcal H}_1  C^{-1/2}$ and $\frac{\kappa_3}{6\sigma ^3}$ are indeed equal. Recall
from the previous section that 
\begin{equation*}
\widetilde{\mathcal H}_1 = \frac{\int_0^1 \Xi_s^{14} ds}{2\int_0^1 \Xi_s^{11} ds}, \qquad \int_0^1 \Xi_s^{11} ds = C,
\qquad \int_0^1 \Xi_s^{14} ds = \E 
\left[ f( W_1 )\int_{0}^{1} \E^2[f'(W_1)| \mathcal F_s] ds \right].  
\end{equation*} 
Hence, we just need to prove the identity
\begin{equation*}
\kappa_3 =3 \E \left[ f( W_1 )\int_{0}^{1} \E^2[f'(W_1)| \mathcal F_s] ds \right].
\end{equation*} 
But $\kappa_3 = \E[f^3(W_1)]$ and It\^o formula  implies that 
\begin{eqnarray*}
3 \E \left[ f( W_1 )\int_{0}^{1} \E^2[f'(W_1)| \mathcal F_s] ds \right] &=& 
3 \E \left[ \int_{0}^{1} \E[f(W_1)| \mathcal F_s] \E^2[f'(W_1)| \mathcal F_s] ds \right] \\[1.5 ex]
&=& \E[f^3(W_1)]
\end{eqnarray*} 
due to the identity $f(W_1)= \int_{0}^{1} \E[f'(W_1)| \mathcal F_s] dW_s$.   \qed  
\end{ex}

\section{Proofs} \label{proof}
\setcounter{equation}{0}
\subsection{A stochastic expansion}
Below, we denote by $K$ a generic positive constant, which may change from line to line. We also write $K_p$ if the constant 
depends on an external parameter $p$. \newline \newline 
{\it Proof of Theorem \ref{th3}:} 
First, we remark that all processes of the type $(b^{[k_1 \ldots k_m]}_s)_{s\geq 0}$ ($k_j\in \{1,2\}$), which we consider
below, are continuous and so locally bounded. Applying the localization technique described in 
Section 3 of \cite{BGJPS}
we can assume w.l.o.g. that these processes are bounded in $(\omega,s)$, which we do from now on.   
We decompose 
\begin{equation*}
\De_n^{-1/2} \Big(V_n(f) - V(f) \Big) = M_n + R_n^{(1)} + R_n^{(2)}
\end{equation*}
with 
\begin{eqnarray*}
R_n^{(1)} &=& \Delta_n^{-1/2} \Big(V_n(f) -  \Delta_n \sum_{i=1}^{1/\De_n} f (\alpha_i^n)  \Big), \\[1.5 ex]
R_n^{(2)} &=& \Delta_n^{-1/2} \Big(\Delta_n \sum_{i=1}^{1/\De_n} \rho_{b^{[1]}_{t_{i-1}}} - V(f)  \Big).
\end{eqnarray*}
We start with the asymptotic expansion of the quantity $R_n^{(2)}$. Due to Burkholder inequality
any process $Y$ of the form (\ref{sde}) with bounded coefficients $b^{[1]}$, $b^{[2]}$  satisfies the inequality
\begin{equation} \label{burkholder}
\E[|Y_t - Y_s|^p]\leq C_p |t-s|^{p/2}
\end{equation}
for any $p\geq 0$. In particular, this inequality holds for the processes $b^{[1]}$, $b^{[2]}$, $b^{[2.2]}$,
$b^{[2.1]}$, $b^{[1.2]}$, $b^{[1.1]}$ as they are diffusion processes (due to It\^o formula). Applying (\ref{burkholder})
and the Taylor expansion we deduce that 
\begin{eqnarray*}
R_n^{(2)} &=& \Delta_n^{-1/2} \sum_{i=1}^{1/\De_n} \int_{t_{i-1}}^{t_{i}}
\{\rho_{b^{[1]}_{t_{i-1}}} - \rho_{b^{[1]}_{s}} \} ds \\[1.5 ex] 
&=& \Delta_n^{-1/2} \sum_{i=1}^{1/\De_n} \int_{t_{i-1}}^{t_{i}}
\{\rho'_{b^{[1]}_{t_{i-1}}} (b^{[1]}_{t_{i-1}} - b^{[1]}_s)  - \frac 12
\rho''_{b^{[1]}_{t_{i-1}}} (b^{[1]}_{t_{i-1}} - b^{[1]}_s)^2\} ds \\[1.5 ex]
&+& o_{\mathbb P}(\Delta_n^{1/2}) \\[1.5 ex]
&=:& R_n^{(2.1)} + R_n^{(2.2)} + o_{\mathbb P}(\Delta_n^{1/2}).  
\end{eqnarray*}
Recall that 
\begin{equation*}
b^{[1]}_t = b^{[1]}_0 + \int_0^t b^{[1.2]}_sds + \int_0^t b^{[1.1]}_sdW_s.
\end{equation*}
We conclude the identity
\begin{eqnarray*}
R_n^{(2.2)} &=&   -\frac{\Delta_n^{-1/2}}{2} \sum_{i=1}^{1/\De_n} \rho''_{b^{[1]}_{t_{i-1}}} \int_{t_{i-1}}^{t_{i}}
(b^{[1]}_{t_{i-1}} - b^{[1]}_s)^2 ds \\[1.5 ex]
&=&    -\frac{\Delta_n^{-1/2}}{2} \sum_{i=1}^{1/\De_n} \rho''_{b^{[1]}_{t_{i-1}}} 
|b^{[1.1]}_{t_{i-1}}|^2 \int_{t_{i-1}}^{t_{i}}
(W_{t_{i-1}} - W_s)^2 ds  + o_{\mathbb P}(\Delta_n^{1/2}) \\[1.5 ex]
&=:&  \Delta_n^{1/2} (N_{n,4}^{(1)} + o_{\mathbb P}(1)).
\end{eqnarray*}
For the term $R_n^{(2.1)}$ we obtain the decomposition
\begin{eqnarray*}
R_n^{(2.1)} &=& -\Delta_n^{-1/2} \sum_{i=1}^{1/\De_n} \rho'_{b^{[1]}_{t_{i-1}}}  \int_{t_{i-1}}^{t_{i}}
\Big( \int_{t_{i-1}}^{s} b^{[1.2]}_u du +   \int_{t_{i-1}}^{s} b^{[1.1]}_u dW_u\Big) ds \\[1.5 ex]
&=&  -\Delta_n^{-1/2} \sum_{i=1}^{1/\De_n}  \rho'_{b^{[1]}_{t_{i-1}}} \Big( \frac{\Delta_n^2}{2} 
b^{[1.2]}_{t_{i-1}} + b^{[1.1]}_{t_{i-1}} \int_{t_{i-1}}^{t_{i}} (W_s - W_{t_{i-1}}) ds\Big) \\[1.5 ex]
&+& o_{\mathbb P}(\Delta_n^{1/2}) \\[1.5 ex]
&=:& \Delta_n^{1/2} (N_{n,4}^{(2)} + N_{n,5} + o_{\mathbb P}(1)).   
\end{eqnarray*} 
We remark that 
\begin{equation*}
N_{n,4} = N_{n,4}^{(1)}+ N_{n,4}^{(2)}.
\end{equation*}
The treatment of the quantity $R_n^{(1)}$ is a bit more involved. We apply again (\ref{burkholder}) and Taylor expansion:
\begin{eqnarray*}
R_n^{(1)} &=& \Delta_n^{1/2} \sum_{i=1}^{1/\De_n} \left(
f'(\alpha_i^n)\{\frac{\De_i^n X}{\sqrt{\De_n}} - \alpha_i^n\} + \frac 12 f''(\alpha_i^n)
\{\frac{\De_i^n X}{\sqrt{\De_n}} - \alpha_i^n\}^2 \right) + o_{\mathbb P}(\Delta_n^{1/2}) \\[1.5 ex]
&=:& R_n^{(1.1)} + R_n^{(1.2)} + o_{\mathbb P}(\Delta_n^{1/2}). 
\end{eqnarray*} 
For the term $R_n^{(1.2)}$ we obtain the decomposition
\begin{eqnarray*}
R_n^{(1.2)} &=& \frac{\Delta_n^{-1/2}}{2} \sum_{i=1}^{1/\De_n}
 f''(\alpha_i^n) \Big( \int_{t_{i-1}}^{t_{i}} b^{[2]}_s ds + \int_{t_{i-1}}^{t_{i}} 
b^{[1]}_s - b^{[1]}_{t_{i-1}} dW_s\Big)^2 \\[1.5 ex]
&=& \frac{\Delta_n^{-1/2}}{2} \sum_{i=1}^{1/\De_n}
 f''(\alpha_i^n) \Big( \Delta_n b^{[2]}_{t_{i-1}} + b^{[1.1]}_{t_{i-1}} \int_{t_{i-1}}^{t_{i}} 
(W_s - W_{t_{i-1}}) dW_s\Big)^2 + o_{\mathbb P}(\Delta_n^{1/2}) \\[1.5 ex]
&=& \frac{\Delta_n^{3/2}}{2} \sum_{i=1}^{1/\De_n}
 f''(\alpha_i^n) \Big(  b^{[2]}_{t_{i-1}} + \frac 12 b^{[1.1]}_{t_{i-1}} 
H_2 \Big(\frac{\Delta_i^n W}{\sqrt{\Delta_n}} \Big)\Big)^2 + o_{\mathbb P}(\Delta_n^{1/2}) \\[1.5 ex]
&=& \Delta_n^{1/2} (N_{n,3} + o_{\mathbb P}(1)).
\end{eqnarray*} 
The quantity $R_n^{(1.1)}$ is decomposed as 
\begin{eqnarray*}
R_n^{(1.1)} &=& \sum_{i=1}^{1/\De_n}
 f'(\alpha_i^n) \Big( \int_{t_{i-1}}^{t_{i}} b^{[2]}_s ds + \int_{t_{i-1}}^{t_{i}} 
\{b^{[1]}_s - b^{[1]}_{t_{i-1}} \} dW_s\Big) \\[1.5 ex] 
&=&  R_n^{(1.1.1)} + R_n^{(1.1.2)}
\end{eqnarray*}
with 
\begin{eqnarray*}
R_n^{(1.1.1)} &=& \Delta_n \sum_{i=1}^{1/\De_n}
 f'(\alpha_i^n) \Big( b^{[2]}_{t_{i-1}} ds + \frac 12 b^{[1.1]}_{t_{i-1}} 
H_2 \Big(\frac{\Delta_i^n W}{\sqrt{\Delta_n}} \Big) \Big), \\[1.5 ex] 
R_n^{(1.1.2)} &=& 
\sum_{i=1}^{1/\De_n}
 f'(\alpha_i^n) \Big( \int_{t_{i-1}}^{t_{i}} \{b^{[2]}_s- b^{[2]}_{t_{i-1}} \} ds \\[1.5 ex] 
&+& 
\int_{t_{i-1}}^{t_{i}} \Big( \int_{t_{i-1}}^{s} b^{[1.2]}_u du + 
\int_{t_{i-1}}^{s} \{b^{[1.1]}_u - b^{[1.1]}_{t_{i-1}}  \}dW_u \Big) dW_s\Big).
\end{eqnarray*}
We remark that 
\begin{equation*}
R_n^{(1.1.1)} = \Delta_n^{1/2} N_{n,1}.
\end{equation*}
Since $f'$ is an odd function (because $f$ is even) we deduce that
\begin{eqnarray*}
R_n^{(1.1.2)} &=& 
\sum_{i=1}^{1/\De_n}
 f'(\alpha_i^n) \Big( b^{[2.1]}_{t_{i-1}} \int_{t_{i-1}}^{t_{i}} (W_s- W_{t_{i-1}}) ds 
+ \frac{\Delta_n^{3/2} b^{[1.1.1]}_{t_{i-1}}}{6} H_3 \Big(\frac{\Delta_i^n W}{\sqrt{\Delta_n}} \Big) \\[1.5 ex] 
&+& b^{[1.2]}_{t_{i-1}} \int_{t_{i-1}}^{t_{i}} (s- t_{i-1}) dW_s \Big)
+ o_{\mathbb P}(\Delta_n^{1/2}).
\end{eqnarray*}
As
\begin{equation*}
R_n^{(1.1.2)} = \Delta_n^{1/2} (N_{n,2} + o_{\mathbb P}(1)),
\end{equation*} 
we are done. \qed

\section{Appendix}
\setcounter{equation}{0}

In this subsection we present a law of large numbers and a multivariate functional stable convergence theorem, 
which is frequently used in this paper. 
For any $k=1,\ldots, d$, let $g_k: C([0,1]) \rightarrow \R$ be a measurable function with polynomial growth, i.e.
\beas
|g_k(x)|\leq K(1+ \|x\|_\infty^p),
\eean 
for some $K>0$, $p>0$ and $\|x\|_\infty = \sup_{z\in [0,1]} |x(z)|$. In most cases $g_k$ will be a function of $x(1)$;
the path-dependent version is only required to account for the asymptotic behaviour of the functional $C_n$. 
Let $(a_s)_{s\geq 0}$ be an $\R^d$-valued, 
$(\mathcal F_s)$-adapted, continuous and bounded stochastic process. Our first result is the following theorem.

\begin{theo} \label{generalucp}
Let $g: \R^d \times C([0,1]) \rightarrow \R$ be a measurable function with polynomial growth in the last variable
and $a=(a_1, \ldots, a_d)$. Then
it holds that
\beas
\De_n\sum_{i=1}^{1/\De_n} g\Big(a_{t_{i-1}},  \De_n^{-1/2}
\{W_{t_{i-1} + s\De_n} - W_{t_{i-1}} \}_{0\leq s\leq 1}\Big) \toop \int_0^1 \rho (a_s,g) ds 
\eean 
with
\beas
\rho (z,g) := \E[g(z,  
\{W_s \}_{0\leq s\leq 1})], \qquad z\in \R^d.
\eean 
\end{theo} 
{\it Proof of Theorem \ref{generalucp}:} Since $
\De_n^{-1/2}
\{W_{t_{i-1} + s\De_n} - W_{t_{i-1}} \}_{0\leq s\leq 1} \eqschw \{W_s \}_{0\leq s\leq 1}. 
$ we obtain that 
\beas
&&\De_n \sum_{i=1}^{1/\De_n} \E\Big[g\Big(a_{t_{i-1}},  \De_n^{-1/2}
\{W_{t_{i-1} + s\De_n} - W_{t_{i-1}} \}_{0\leq s\leq 1}\Big)|\mathcal F_{t_{i-1}} \Big] \\[1.5 ex] 
&&= \De_n \sum_{i=1}^{1/\De_n} \rho (a_{t_{i-1}},g) \toop \int_0^1 \rho (a_s,g) ds. 
\eean 
On the other hand, we deduce that 
\beas
&&\De_n\sum_{i=1}^{1/\De_n} g\Big(a_{t_{i-1}},  \De_n^{-1/2}
\{W_{t_{i-1} + s\De_n} - W_{t_{i-1}} \}_{0\leq s\leq 1}\Big) \\[1.5 ex] 
&&- \De_n \sum_{i=1}^{1/\De_n} \E\Big[g\Big(a_{t_{i-1}},  \De_n^{-1/2}
\{W_{t_{i-1} + s\De_n} - W_{t_{i-1}} \}_{0\leq s\leq 1}\Big)|\mathcal F_{t_{i-1}} \Big] \toop 0,
\eean 
because 
\beas
\De_n^2 \sum_{i=1}^{1/\De_n} \E\Big[g^2\Big(a_{t_{i-1}},  \De_n^{-1/2}
\{W_{t_{i-1} + s\De_n} - W_{t_{i-1}} \}_{0\leq s\leq 1}\Big)|\mathcal F_{t_{i-1}} \Big] \toop 0.
\eean 
This completes the proof. \qed \\ \\
Next, we consider a sequence of $d$-dimensional processes
$Y_t^n = (Y_{1,t}^n, \ldots, Y_{d,t}^n)$ defined via
\beas
Y_{k,t}^n &=& \De_n^{1/2} \sum_{i=1}^{[t/\De_n]} a^k_{t_{i-1}} \Big[  g_k\Big( \De_n^{-1/2}
\{W_{t_{i-1} + s\De_n} - W_{t_{i-1}} \}_{0\leq s\leq 1}  \Big)  \\[1.5 ex]
  &-& \E g_k\Big( \De_n^{-1/2}
\{ W_{t_{i-1} + s\De_n} - W_{t_{i-1}} \}_{0\leq s\leq 1}  \Big) \Big], \qquad k=1,\ldots, d. 
\eean 
The stable convergence of $Y^n$ is as follows.

\begin{theo} \label{generalclt}
It holds that 
\beas
Y^n_t \stab Y_t = \int_0^t v_s dW_s + \int_0^t (w_s-v_s v_s^{\star})^{1/2} dW'_s,  
\eean 
where the functional convergence is stable in law, $W'$ is a $d$-dimensional Brownian motion independent of $\mathcal F$,
and the processes $(v_s)_{s\geq 0}$ in $\R^{d}$ and $(w_s)_{s\geq 0}$ in 
$\R^{d \times d}$ are defined as
\beas
v_s^{k} &=& a_s^{k} \E \Big[g_k ( 
\{W_s \}_{0\leq s\leq 1}) W_1 \Big], \\[1.5 ex]
w_s^{kl} &=& a_s^{k} a_s^{l} \mbox{cov}\Big[ g_k ( 
\{W_s \}_{0\leq s\leq 1}),  g_l ( 
\{W_s \}_{0\leq s\leq 1})\Big],
\eean
with $1\leq k,l\leq d$. In particular, it holds that  $\int_0^t w_s^{1/2} dW'_s \sim MN\left(0, \int_0^t w_s ds \right)$.
\end{theo}   
{\it Proof of Theorem \ref{generalclt}:} We write $Y^n_t = \sum_{i=1}^{[t/\De_n]} \chi_i^n$ with 
\beas
\chi_{i,k}^n &=& \De_n^{1/2}  a^k_{t_{i-1}} \Big[  g_k\Big( \De_n^{-1/2}
\{W_{t_{i-1} + s\De_n} - W_{t_{i-1}} \}_{0\leq s\leq 1}  \Big)  \\[1.5 ex]
  &-& \E g_k\Big( \De_n^{-1/2}
\{ W_{t_{i-1} + s\De_n} - W_{t_{i-1}} \}_{0\leq s\leq 1}  \Big) \Big], \qquad k=1,\ldots, d. 
\eean 
According to Theorem IX.7.28 of \cite{JS} we need to show that 
\begin{eqnarray}
&&\label{p1} \sum_{i=1}^{[t/\De_n]} \E[\chi_{i,k}^n \chi_{i,l}^n | \mathcal F_{t_{i-1}}] \toop \int_0^t w_s^{kl} ds, \\[1.5 ex]
&& \label{p2} \sum_{i=1}^{[t/\De_n]} \E[\chi_{i,k}^n \De_i^n W | \mathcal F_{t_{i-1}}] \toop \int_0^t v_s^{k} ds, \\[1.5 ex]
&&\label{p3} \sum_{i=1}^{[t/\De_n]} \E[|\chi_{i,k}^n|^2 1_{\{|\chi_{i,k}^n|>\epsilon \}} | \mathcal F_{t_{i-1}}] 
\toop 0 \qquad \forall \epsilon>0,  \\[1.5 ex]
&&\label{p4} \sum_{i=1}^{[t/\De_n]} \E[\chi_{i,k}^n \De_i^n Q | \mathcal F_{t_{i-1}}] \toop 0,
\end{eqnarray}
where $1\leq k,l\leq d$ and the last condition must hold for all bounded continuous martingales $Q$ with $[W,Q]=0$. Conditions 
(\ref{p1}) and (\ref{p2}) are obvious since
\beas
\De_n^{-1/2}
\{W_{t_{i-1} + s\De_n} - W_{t_{i-1}} \}_{0\leq s\leq 1} \eqschw \{W_s \}_{0\leq s\leq 1}. 
\eean
Condition (\ref{p3}) follows from
\beas
\sum_{i=1}^{[t/\De_n]} \E[|\chi_{i,k}^n|^2 1_{\{|\chi_{i,k}^n|>\epsilon \}} | \mathcal F_{t_{i-1}}]\leq 
\epsilon^{-2} \sum_{i=1}^{[t/\De_n]} \E[|\chi_{i,k}^n|^4 | \mathcal F_{t_{i-1}}]\leq K \De_n\rightarrow 0,
\eean
which holds since the process $a$ is bounded 
and $g_k$ is of polynomial growth. In order to prove the last condition,
we use the It\^o-Clark representation theorem
\beas
  &&g_k\Big( \De_n^{-1/2}
\{W_{t_{i-1} + s\De_n} - W_{t_{i-1}} \}_{0\leq s\leq 1}  \Big)  \\[1.5 ex]
  &&- \E g_k\Big( \De_n^{-1/2}
\{ W_{t_{i-1} + s\De_n} - W_{t_{i-1}} \}_{0\leq s\leq 1}  \Big)  = \int_{t_{i-1}}^{t_{i}} \eta_{k,s}^n dW_s 
\eean 
for some predictable process $\eta_{k}^n$. It\^o isometry implies the identity  
\beas 
\E[\chi_{i,k}^n \De_i^n Q | \mathcal F_{t_{i-1}}] = \De_n^{1/2}  a^k_{t_{i-1}} 
\E\left[\int_{t_{i-1}}^{t_{i}} \eta_{k,s}^n d[W,Q]_s| \mathcal F_{t_{i-1}} \right]=0.
\eean
This completes the proof of Theorem \ref{generalclt}. \qed

\subsection*{Acknowledgment}
Mark Podolskij acknowledges financial support from CREATES funded by the Danish
National Research Foundation. Yoshida's research was in part supported by Grants-in-Aid for 
Scientific Research No. 19340021, No. 24340015 (Scientific Research), No.24650148 (Challenging Exploratory Research); 
the Global COE program "The Research and Training Center for New Development in Mathematics" of the Graduate School of Mathematical Sciences, University of Tokyo; Cooperative Research Program of the Institute of Statistical Mathematics; and by NS Solutions Corporation.
Both authors would like to thank Bezirgen Veliyev for helpful comments on the earlier version of this paper.

\end{document}